\documentclass[pdflatex,sn-mathphys-num]{sn-jnl}

\usepackage{hyperref}%
\usepackage[T1]{fontenc}
\usepackage[utf8]{inputenc}
\usepackage[english]{babel}

\usepackage{graphicx} 
\usepackage{amsmath}
\usepackage{mathrsfs}
\usepackage{amssymb}
\usepackage{mathtools}
\usepackage{amsthm}

\usepackage{booktabs}
\usepackage{tabularx}
\usepackage{verbatim}
\usepackage{multirow}
\usepackage{graphicx} 

\usepackage{textgreek}
\usepackage{xfrac}    
\usepackage{faktor}
\usepackage{xcolor}
\usepackage[normalem]{ulem}
\usepackage{todonotes}

\newcommand{\inner}[2]{\ifthenelse{\equal{#2}{}}{\left\langle\cdot,\cdot\right\rangle_{#1}}{\left\langle#2\right\rangle_{#1}}}
\newcommand{\norm}[2]{\ifthenelse{\equal{#2}{}}{\left\|\cdot\right\|_{#1}}{\left\|#2\right\|_{#1}}}
\newcommand{\seminorm}[2]{\ifthenelse{\equal{#2}{}}{\left|\cdot\right|_{#1}}{\left|#2\right|_{#1}}}
\newcommand{\calh}{\mathcal H}
\newcommand{\calx}{\mathcal X}
\newcommand{\N}{\mathbb{N}}
\newcommand{\R}{\mathbb{R}}
\newcommand{\Domain}{D}
\newcommand{\Range}{R}
\newcommand{\Null}{N}
\newcommand{\fa}{\hbox{ for all }}
\DeclareMathOperator{\diag}{diag}
\DeclareMathOperator{\Sp}{span}
\DeclareMathOperator*{\argmin}{argmin}

\newtheorem{definition}{Definition}
\newtheorem{proposition}{Proposition}
\newtheorem{corollary}{Corollary}
\newtheorem{lemma}{Lemma}

\newtheorem{theorem}{Theorem}
\newtheorem{remark}{Remark}
\newtheorem{assumption}{Assumption}

\begin{document}

\title[]{Weak convergence rates for spectral regularization via sampling inequalities}


\author*[1]{\fnm{Sabrina} \sur{Guastavino}}\email{sabrina.guastavino@unige.it}

\author[2]{\fnm{Gabriele} \sur{Santin}}\email{gabriele.santin@unive.it}

\author[3]{\fnm{Francesco} \sur{Marchetti}}\email{francesco.marchetti@unipd.it}
\author[1]{\fnm{Federico} \sur{Benvenuto}}\email{federico.benvenuto@unige.it}

\affil*[1]{\orgdiv{Department of Mathematics}, \orgname{University of Genova}, \orgaddress{\city{Genova}, \country{Italy}}}

\affil[2]{\orgdiv{Department of Environmental Sciences, Informatics and Statistics}, \orgname{Ca’ Foscari University of Venice}, \orgaddress{\city{Venice}, \country{Italy}}}

\affil[3]{\orgdiv{Department of Mathematics ``Tullio Levi-Civita''}, \orgname{University of Padova}, \orgaddress{\city{Padova}, \country{Italy}}}




\abstract{
Convergence rates in spectral regularization methods quantify the approximation error in inverse problems as a function of the noise level or the number of sampling points.
Classical strong convergence rate results typically rely on source conditions, which are essential for estimating the truncation error.
However, in the framework of kernel approximation, the truncation error in the case of Tikhonov regularization can be characterized entirely through sampling inequalities, without invoking source conditions.
In this paper, we first generalize sampling inequalities to spectral regularization, 
and then, by exploiting the connection between inverse problems and kernel approximation, we derive weak convergence rate bounds for inverse problems, independently of source conditions. 
These weak convergence rates are established and analyzed when the forward operator is compact and uniformly bounded, or the kernel operator is of trace class.}




\maketitle

\section{Introduction}

The problem of recovering a function from discrete data arises in various application domains such as inverse problems, interpolation, kernel approximation and learning \cite{wahba1977practical,engl1996regularization, bertero1985linear,vapnik2013nature,cucker_learning_2007}.
The quality of the estimation is typically assessed in terms of convergence rates, which quantify the decay of the approximation error as the number of data samples increases.
In inverse problems, given a finite and discrete set of noise-corrupted observations, spectral regularization schemes yield one-parameter families of estimators whose squared-norm error satisfies the classical bias–variance trade-off.
The key point in the study of convergence rates is the \emph{source condition}, which is crucial for bounding the bias term.
Source conditions encode a priori smoothness assumptions on the (possibly generalized) solution \( f^\dagger \) of an inverse problem described by the forward equation \( Af^\dagger = g \), where $A: \mathcal{H}_1 \to \mathcal{H}$ is a compact linear operator.
The classical form is
\begin{equation*}
f^\dagger = \psi(A^*A)\,w, \quad w \in \mathcal{H}_1,
\end{equation*}
where $\psi$ is an index function, for example $\psi(t) = t^\mu$ for the classical Holder-type source conditions
\cite{engl1996regularization,loubes2009review}.
In inverse problems, convergence rate bounds have been established in both a deterministic setting, where the noise is infinite dimensional, and a statistical setting, where the data are drawn from a (possibly unknown) distribution.
These analyses accommodate general classes of operators, such as Hilbert–Schmidt or trace-class operators \cite{wang_optimal_2011,blanchard2017optimal}.
Moreover, a precise correspondence between the convergence rates in these two settings has been established \cite{guastavino2020convergence}.
The main limitation of source conditions is that, while they guarantee the theoretical  computability of the bounds, they impose assumptions that are often unrealistic and rarely satisfied in practice.

A different but related line of research arises in kernel approximation and approximation theory, where convergence rate results are derived under a deterministic setting, and sampling is therefore treated as a geometric process yielding a set of points $\mathcal{X}_n = \{x_1,\ldots,x_n\} \subset \mathcal X$ contained in the compact subset over which functions are defined.
The central tool for bounding the approximation error are the \textit{sampling inequalities} that are defined in terms of the \emph{fill distance} 
\[
h_{\mathcal{X}_n,\mathcal{X}} := \sup_{x \in \mathcal{X}} \min_{x_i \in \mathcal{X}_n} \|x - x_i\|~.
\]
This quantity reflects the density of the distribution of the sampling points, and offer a concrete way to measure how well the domain is covered \cite{Rieger10,Rieger17,Rieger20}.
Since the deterministic approximation error studied in approximation theory mirrors the bias term in statistical inverse problems, it is natural to ask:
\begin{center}
\textit{Is it possible to establish error estimates for spectral regularization of inverse problems by using sampling inequalities and, therefore, independently of source conditions?} 
\end{center}
In this study, we provide a positive response to this question under the condition that the error estimates are given in a weak sense.
To achieve this result, we analyze the interplay between operator-theoretic regularization and kernel approximation theory by exploiting that the image of a compact operator is a Reproducing Kernel Hilbert Space (RKHS). 
Our main contributions are the following:
\begin{itemize}
    \item[1.] we generalize the bounds for the approximation error using sampling inequalities to the setting of \emph{spectral regularization methods}, obtaining deterministic bias bounds depending explicitly on the fill distance \(h\);
    \item[2.]  we derive upper bounds on the variance term under gaussian noise models, by incorporating additive noise in the sampling process;
    \item[3.] we use such approximation bounds to give weak error estimates for spectral regularized solutions in inverse problems.
\end{itemize} 

Weak error estimates are obtained by selecting suitable test functions.
In this study we consider two classes of test functions, each of which depends on the forward operator, providing two distinct ways to evaluate the convergence of spectral reconstruction error in a weak sense.
The main advantage of using sampling inequalities to derive these error bounds is that the resulting estimates do not require any {\it a priori} assumptions on the smoothness of the source, which is rarely available in practical applications.
In contrast, such weak estimates depend on geometric properties of the sampling set, most notably the fill distance.
This shift can make the resulting error bounds potentially meaningful in real-world problems, since the fill distance can either be computed directly from the data, as in the case of a training set in supervised learning, or be known in advance, as in image reconstruction where the sampling pattern is controlled (e.g., sensor locations, acquisition grids).

The paper is organized as follows. In Section \ref{sec:mathematical_setup} we present the mathematical setup and highlight the connection between approximation problems in RKHS and inverse problems solutions.
In Section \ref{sec:main_results}, we state our main results, which include new bounds for spectral regularized solutions of approximation problems, obtained via sampling inequalities, and their translation into weak error estimates for spectral regularized solutions in inverse problems.
Section \ref{sec:applications} provides application examples. In Sections \ref{sec:proof_approx_err} and \ref{sec:proof_noise} the proofs of the main results are presented, and finally the conclusions are drawn in Section \ref{sec:cocnlusions}.

\section{Mathematical setup}\label{sec:mathematical_setup}

\subsection{Compact linear operators and RKHSs}

In this section, we recall some basic facts about linear compact operators and their connection with RKHSs, with the aim of fixing the notation used throughout the paper. We refer to \cite{engl1996regularization} and~\cite{Desoer1963} for a detailed treatment of these topics.

\subsubsection{Compact linear operators}

Let us consider a compact linear operator $A:\mathcal H_1 \to \mathcal H_2$, where $\mathcal{H}_1$ is a Hilbert space and $\mathcal H_2$ is a Hilbert space of functions over a set $\mathcal X$, and denote by $\Domain(A)$, $\Range(A)$ and $\Null(A)$ the domain, range, and null space of $A$, respectively. 
In this setting there exists a singular system $\{(\sigma_j, u_j,v_j)\}_j$ of $A$, i.e., a $\calh_1$-orthonormal system $\{v_j\}_{j\in\N}\subset \calh_1$, a $\calh_2$-orthonormal system $\{u_j\}_{j\in\N}\subset 
\calh_2$, and strictly positive and non increasing numbers $\{\sigma_j\}_{j\in\N}$ such that $A$ has the spectral decomposition
\begin{equation*}
A v = \sum_{j=1}^\infty \sigma_j \inner{\calh_1}{v_j, v} u_j, \;\; v\in \calh_1.
\end{equation*}
It follows that the adjoint $A^*: \calh_2\to\calh_1$ of $A$ is given by
\begin{equation*}
A^* u = \sum_{j=1}^\infty \sigma_j \inner{\calh_2}{u_j, u} v_j, \;\; u\in \calh_2,
\end{equation*}
so that
\begin{equation}\label{eq:range_a_adj}
\Range(A^*)=\left\{\sum_j a_j v_j\in \calh_1 : \sum_j v_j^2/\sigma_j^2<+\infty\right\}.
\end{equation}
Moreover, provided that the expansion
\begin{equation}\label{eq:phi_x}
\phi_x \coloneqq \sum_{j=1}^{\infty} \sigma_j u_j(x) v_j
\end{equation}
converges for each $x \in \calx$, using the fact that $\{v_j\}_j$ is a basis of $\Null(A)^\perp$ the operator $A$ can be written as
\begin{equation}
\label{eq:A_via_phi}
(Aw)(x) = \langle  w, \phi_x\rangle_{\mathcal{H}_1},\;w\in \calh_1.
\end{equation}
We remark that $\phi_x\in \Null(A)^\perp$ and, letting $\phi : \mathcal X \longrightarrow \mathcal H_1$ 
be a feature map defined by $\phi(x) = \phi_x$ for any $x\in\calx$, it holds that 
\begin{equation}\label{eq:feature_space}
\mathcal{F}\coloneqq\overline{\textrm{span}\{\phi(x) ~, ~ x\in\mathcal{X}\}} = N(A)^{\perp}.
\end{equation}

We recall that the pseudoinverse is a bounded operator if and only if $A$ has closed range and, since $A$ is compact, this is true if and only if it has a 
finite dimensional range.
The pseudo inverse of $A$ has the representation
\begin{equation*}
A^\dagger u = \sum_{j=1}^\infty \frac1{\sigma_j} \inner{\calh_2}{u_j, u} v_j, \;\; u\in \Range(A) \otimes \Range(A)^\perp,
\end{equation*}
and its adjoint, which satisfies $(A^\dagger)^* = (A^*)^\dagger$, is
\begin{equation}\label{eq:A_adj_ps}
(A^\dagger)^* v = \sum_{j=1}^\infty \frac1{\sigma_j} \inner{\calh_1}{v_j, v} u_j, \;\; v\in \Range(A^\dagger) \otimes \Range(A^\dagger)^\perp = \Null(A) \otimes \Null(A)^\perp.
\end{equation}

\subsubsection{Spectral regularization}

Spectral regularization tools play a fundamental role in inverse problems.
In general, a spectral regularization is a map $\mathfrak{R} : \mathcal B \times \mathcal H_2 \times \mathbb{R}_{+} \to \mathcal H_1$ defined by
\begin{equation}
\label{spectral reg sol inv pb}
    \mathfrak{R}(B,z,\lambda) := s_{\lambda}(B^*B)B^*z,
\end{equation}
where $B\in\mathcal{B}$ is a bounded linear operator, $z\in\mathcal{H}_2$, $\lambda\in\mathbb{R}_+$, $s_{\lambda}:[0,\infty)\to\mathbb{R}$ denotes the {\it regularization function} and the operator $s_{\lambda}(B^*B)$ is defined via the spectral calculus. 
In order to ensure stability and convergence of the spectral regularization, the family $\{s_{\lambda}\}_{\lambda>0}$ has to satisfy the following:
\begin{assumption}\label{def:regularization_fun}
The regularization (or filtering) function $s_{\lambda}:[0,\infty)\to \mathbb{R}$ for $\lambda>0$ 
satisfies the 
following properties:
\begin{enumerate}
\item there exists a constant $D>0$ such that
\begin{equation}\label{property1 s_lambda}
\sup_{t\in [0,\Vert B\Vert^2]} |t s_{\lambda}(t)|\le D ~ ~ \text{ uniformly in } ~ \lambda>0 ~ ,
\end{equation}
\item there exists a constant $E>0$ such that\
\begin{equation}\label{property2 s_lambda}
\sup_{\lambda>0} \sup_{t\in [0,\Vert B\Vert^2]}|\lambda s_{\lambda}(t)|\le E ~ ,
\end{equation}
\item there exists $q>0$ called qualification of the method and constants $C_{a}>0$ such that
\begin{equation}\label{qualification reg function}
\sup_{t\in [0,\Vert B\Vert^2]} |t^{a}(1-t s_{\lambda}(t))| \le C_{a} \lambda^{a} ~ ~ \forall ~ \lambda>0 ~ \text{ and } ~ ~ 0\le a \le q.
\end{equation}
\end{enumerate}
\end{assumption}

The idea of spectral regularization is to provide approximated solutions of a linear operator equation with noisy data by filtering out high frequencies.
Typical examples are the following.
\begin{itemize}
    \item Tikhonov regularization: in this case the regularization function is given by $s_{\lambda}(t)=(\lambda + t)^{-1}$ and the qualification is $q=1$.
    \item Truncated singular value decomposition (or spectral cut-off): in this case the regularization function is given by 
    \begin{equation}
        s_{\lambda}(t)=\begin{cases}
        \frac{1}{t} ~, & \text{ if } t\ge\lambda \\
        0 ~, & \text{ if } t<\lambda
        \end{cases}
    \end{equation}
    and $q$ is arbitrary.
    \item Landweber regularization: given $\gamma$ such that $|1-\gamma t|<1$ for $t\in\sigma(B^*B)$, the regularization function is given by 
    \begin{equation}
        s_{\lambda}(t) = t^{-1} \left( 1-(1-\gamma t )^{(1/\lambda)} \right) 
    \end{equation}
    and $q$ is arbitrary.\end{itemize} 

\subsubsection{Linear compact operators and RKHSs}

A key ingredient in our results is the connection between linear compact operators and RKHSs.
We first recall the definition of RKHS.
\begin{definition}\label{def:rkhs}
Let $\calh_K\coloneqq\calh_K(\calx)$ be an Hilbert space of real valued functions on a non-empty set $\mathcal{X}$. 
$\calh_K$ is said to be a RKHS if there exists a reproducing kernel $K:\mathcal{X}\times\mathcal{X}\to\R$, i.e., $K(x,\cdot) \in \mathcal H$ and 
$g(x)=\langle g,K(\cdot, x)\rangle_{\calh_K}$.
\end{definition}

Equivalently, $\calh_K$ is a RKHS if the evaluation functionals $L_x:g\in\mathcal{H}\to L_x(g):=g(x)$ are continuous for all $x\in\mathcal{X}$, and in this case $K(\cdot, x)$ is the Riesz representer of $L_x$ for all $x\in\calx$.
The space $\calh_K$ is often named the native space of $K$ in the literature on approximation theory~\cite{Wendland05,Fasshauer07,Fasshauer15}.
The following fact is a well-established result in the theory of linear operators (see Chapter 17.3 \cite{kress1999linear}). 
For the sake of clarity, we give a simple proof of this result in our framework.

\begin{theorem}
Let $A$ be a compact linear operator defined by equation \eqref{eq:A_via_phi}. The range $\Range(A)$ of $A$ is a RKHS 
$\calh_K$ with
inner product
\begin{equation}\label{eq:inner_rkhs}
\inner{\calh_K}{g, h}\coloneqq \inner{\calh_1}{A^\dagger g, A^\dagger h},\;\; g,h\in \calh_K,   
\end{equation}
and 
kernel 
\begin{equation}\label{eq:kernel}
K(x,x')\coloneqq \inner{\mathcal{H}_1}{\phi_x,\phi_{x'}},\;\; x,x'\in\calx,
\end{equation}
where $\phi_x$ is defined in equation \eqref{eq:phi_x}.
\end{theorem}
\begin{proof}
First, the expression \eqref{eq:inner_rkhs} is well defined since $\calh_K=\Range(A)\subset\Domain(A^\dagger)$, and it indeed defines an inner product since it 
is clearly symmetric, linear, and positive semi-definite.
The space $\left(\calh_K, \inner{\calh_K}{}\right)$ is thus a pre-Hilbert space. It is in fact complete, since any Cauchy sequence $\{g_n\}_n\subset \calh_K$ 
defines a Cauchy sequence $\{A^\dagger g_n\}_n\subset \calh_1$. Since $\calh_1$ is complete, this sequence has in turn a limit $g\in \calh_1$, and it is 
easy to see 
that $g_n\to A g$ in $\calh_K$.

To see that this space is a RKHS we show that $L_x$ is bounded for all $x\in\calx$ (see Definition~\ref{def:rkhs}). Indeed, for all $u\in \calh_K$ we have 
that there is $v\in \Null(A)^\perp \subset\calh_1$ with $Av=u$, and in particular for all $x\in\calx$ it holds
\begin{align*}
|u(x)| &= |Av(x)| = \left|\inner{\calh_1}{\phi_x, v}\right|\leq \norm{\calh_1}{v}\norm{\calh_1}{\phi_x},\\
\norm{\calh_K}{u}&=\norm{\calh_1}{A^\dagger u} = \norm{\calh_1}{v},
\end{align*}
and thus $\norm{\calh_1'}{L_x}\leq \norm{\calh_1}{\phi_x}<\infty$.

Finally, the definition~\eqref{eq:kernel} and~\eqref{eq:A_via_phi} show that for all $x,x'\in\calx$ we have $K(x, x')=\inner{\mathcal{H}_1}{\phi_x,\phi_{x'}} = 
A \phi_{x'}(x)$, and thus $K(\cdot, x)\in \calh_K$ since $\calh_K=\Range(A)$. 
Moreover, for $u = Av\in \calh_K$ with $v\in \Null(A)^\perp$ we can write
\begin{align*}
\inner{\calh_K}{K(\cdot, x), u}
&= \inner{\calh_K}{A \phi_x, u}
= \inner{\calh_1}{A^\dagger A \phi_x, A^\dagger u}
= \inner{\calh_1}{\phi_x, v}\\
&= A v(x) = u(x),
\end{align*}
which shows that $K$ is the reproducing kernel of $\calh_K$.
\end{proof}
Observe that the definition~\eqref{eq:kernel} also gives the expression
\begin{equation}
\label{eq:svd-kernel}
K(x,x') = \sum_{j=1}^{\infty} \sigma_j^2 u_j(x) u_j(x'),
\end{equation}
and that the theorem shows the existence of an isomorphism between the space $\mathcal F$ in \eqref{eq:feature_space} and the RKHS. 
Namely, the first isomorphism theorem applied to the map
\begin{equation*}
A : \mathcal H_1 \longrightarrow \Range(A) = \mathcal H_K \subseteq \mathcal H
\end{equation*}
yields
\begin{equation*}
\mathcal{F} = N(A)^{\perp} \simeq \faktor{\mathcal H_1}{N(A)} \simeq \Range(A) = \mathcal H_K ~.
\end{equation*}
As the kernel is invariant under unitary transformations, the map $\phi_x$ associated with a given kernel $K$ is determined only up to the action of the unitary group on $\mathcal{H}_1$, i.e.
\begin{equation*}
\label{K}
K(x,x')=\langle \phi_x,\phi_{x'}\rangle_{\mathcal{H}_1}=\langle U\phi_x,U\phi_{x'}\rangle_{\mathcal{H}_1},
\end{equation*}
for each unitary operator $U$ acting on $\mathcal{H}_1$. 
%

\subsection{Semi-discrete inverse problems in RKHSs}\label{sec:kernels_in_inv_prob}

Given any set of nodes $\mathcal{X}_n = \{x_1,\dots, x_n\}\subseteq \mathcal{X}$ and a function $f\in\mathcal{H}_1$, we consider the \textit{semi-discrete operator} operator 
\begin{equation}
(A_n f)_i \coloneqq \langle f,\phi_{x_i} \rangle_{\mathcal{H}_1} = (Af)(x_i),
\end{equation}
which maps $A_n:\mathcal{H}_1 \to \mathbb{R}^n$.
With the normalized inner product on $\mathbb R^n$, its adjoint is
\begin{equation}
    A^*_n  w = \frac{1}{n}\sum_{i=1}^n w_i \phi_{x_i} ~~ \text{ and } ~~
    A_n^*A_n = \frac{1}{n}\sum_{i=1}^n \langle \cdot,\phi_{x_i}\rangle_{\mathcal{H}_1}  \phi_{x_i},
\end{equation}
for $w\in\mathbb{R}^n$.
In this framework, we start by considering the reconstruction of the function
$$
g\coloneqq Af: \mathcal X \to \mathbb R
$$
from its sampled values $y=(y_1,\ldots,y_n)$ with $y_i=g(x_i)$.
In this case, the kernel-based interpolant takes the form
\begin{equation}\label{eq:interpolant}
\hat{g}^{(n)} = k(\mathcal{X}_n)\cdot c ~~\text{where}~~     k(\mathcal{X}_n)\coloneqq (K(\cdot,x_1),\dots,K(\cdot,x_n))
\end{equation}
and the coefficients $c=(c_1,\dots,c_n)$ solves the linear system 
\begin{equation}\label{eq:system}
\mathsf{K}c=y, ~~\text{where}~~ \mathsf{K}=(\mathsf{K}_{ij})=K(x_i,x_j)=(A_n A_n^*)_{ij}
\end{equation}
for $i,j=1,\dots,n$, is the so-called interpolation, or collocation, or simply kernel, or Gram's matrix.
Therefore, the interpolant $\hat{g}^{(n)}$ belongs to $\mathcal{H}_{\mathcal{X}_n}:=\mathrm{span}\{ K(x,\cdot), ~ x\in\mathcal{X}_n\}\subset\mathcal{H}_K$.
When $\mathsf{K}$ is strictly positive definite, the interpolant \eqref{eq:interpolant} is the unique minimizer of the \textit{empirical risk}
\begin{equation}
\label{approx_probl_sampling}
\hat{g}^{(n)} =\argmin_{z\in\mathcal{H}_K} \sum_{i=1}^n 
(y_i-z(x_i))^2 ~.
\end{equation}
We refer to \cite{Fasshauer07,Fasshauer15,Wendland05} for more details on this topic, including classical results.
The associated inverse problem consists in recovering the minimum-norm solution
\begin{equation}
\label{inv_probl_sampling}
\hat{f}^{(n)}:=\argmin_{f'\in F} 
\| f' \|_{\mathcal H_1},
\end{equation}
where
\begin{equation}
\label{f_R}
F \coloneqq \argmin_{f\in\mathcal{H}_1} \sum_{i=1}^n 
(y_i-(A_n f)_i)^2 ~.
\end{equation}
Thanks to the representer theorem (see e.g. \cite[Theorem 1]{Scholkopf01}), the solution admits a finite representation
\begin{equation}\label{direct_representer}
\hat{f}^{(n)} = \phi(\mathcal{X}_n)\cdot c ~~~\text{where}~~~ \phi(\mathcal{X}_n)\coloneqq (\phi_{x_1},\dots,\phi_{x_n})
\end{equation}
with the same coefficients $c$ appearing in the interpolation formula \eqref{eq:interpolant}.
Equivalently, the minimizers \eqref{approx_probl_sampling} and \eqref{inv_probl_sampling} are related by
$$
\hat{f}^{(n)} = A^\dagger \hat{g}^{(n)}.
$$

When the data are perturbed by i.i.d. Gaussian noise,
\begin{equation}
\label{eq:noise}
y^\delta = y + \delta, ~~~ \delta_i \sim G(0,\sigma^2),
\end{equation}
it is natural to consider a regularized reconstruction.
In the direct (kernel) formulation, it is common practice to solve the Tikhonov regularized problem
\begin{equation}
\label{approx_probl_sampling_noisy}
\hat{g}^{(\delta,n)}_{\lambda}:=
\argmin_{z\in\mathcal{H}_K} \sum_{i=1}^n (y^\delta_i-z(x_i))^2
+\lambda \| z\|_{\mathcal{H}_K}^2~,
\end{equation}
with $\lambda>0$ the regularization parameter, whose solution has the finite representation
\begin{equation}
\label{direct_representer_reg_1}
\hat{g}^{(\delta,n)}_{\lambda} = k(\mathcal{X}_n)\cdot c_{\lambda}^\delta, ~~~
(\mathsf{K} + \lambda I) \, c_\lambda^\delta=y^{\delta}.
\end{equation}
The inverse Tikhonov regularized problem 
\begin{equation}
\label{inv_probl_sampling_noisy}
\hat{f}^{(\delta,n)}_{\lambda}:=\arg\min_{f\in\mathcal{H}_1} \sum_{i=1}^n (y^\delta_i - (Af)(x_i))^2 +\lambda \| f\|_{\mathcal{H}_1} ~.
\end{equation}
then satisfies
\begin{equation}
\label{direct_representer_reg}
\hat{f}^{(\delta,n)}_{\lambda} = (A^\dagger \hat{g}^{(\delta,n)}_{\lambda} ) = \phi(\mathcal{X}_n)\cdot c_{\lambda}^\delta.
\end{equation}
Therefore, direct and inverse estimators share the same coefficient vector, with only the basis $\{K(\cdot,x_i)\}$ replaced by $\{\phi_{x_i}\}$ through the relation $f=A^\dagger g$.
In particular, the representer theorem makes the correspondence between the two problems completely explicit, independently of whether regularization is applied.

More generally, a spectral regularization map $s_\lambda$ applied to $A_n^* A_n$ yields the solution
\begin{equation}\label{eq:spectral_reg_f}
    \hat{f}^{(\delta,n)}_{s_{\lambda}} = s_{\lambda}(A_n^*A_n) A_n^*  y^{\delta}.
\end{equation}
For such spectral regularization methods, error estimates are provided only when some source condition is given. See, for example, \cite{blanchard2017optimal} and references therein.
Moreover, spectral regularization provides a generalization of the classical smoothing strategy for interpolation, since it defines an approximation function as
\begin{equation}\label{eq:g_ldn}
\hat{g}_{s_\lambda}^{(\delta,n)} \coloneqq A_{s_\lambda}~y^{\delta}
\end{equation}
where $A_{s_\lambda} : \mathbb R^n \to \Range(A)=\calh_K$, defined by
\begin{equation}\label{eq:A_lambda}
    A_{s_\lambda} \coloneqq A ~ s_\lambda(A_n^* A_n) A_n^*, 
\end{equation}
is the semidiscrete influence operator \cite{vogel2002computational}.
Inspired by this direct-inverse connection, in the next section we present some widely used tools from approximation theory that will be employed to compute weak error estimates for inverse problems, independent of source conditions.

\subsection{Kernel-based approximation and sampling inequalities}\label{sec:sampling}

We assume that $\mathcal{X}\subset\mathbb{R}^d$ with norm $\|\cdot\|_2$, and we define the \textit{fill distance} of the set of points $\mathcal{X}_n$ as follows 
\begin{equation*}
    h_{\mathcal{X}_n,\mathcal{X}} \coloneqq  \sup_{ x \in \mathcal{X}} {\min_{ x_i  \in {\mathcal{X}_n}} \lVert x - x_i \lVert_2 },
\end{equation*}
which measures to what extent the domain is filled by data nodes, and the \textit{separation distance}
\begin{equation}\label{eq:sep_dist}
q_{\mathcal{X}_n} \coloneqq \frac12 \min_{  i \neq j} \lVert {x}_i - x_j \lVert_2,    \end{equation}
which measures the minimal distance between two distinct nodes. 
It is clear that $q_{\calx_n}\leq h_{\calx_n,\calx}$. The reverse inequality characterizes certain sets that are of large interest in kernel-based methods.
\begin{definition}\label{def:quasi_uniform}
Let $\calx\subset\R^d$ be bounded. A sequence $\{\calx_n\}_{n\in\N}$ of sets $\calx_n\subset\calx$, each of $n$ pairwise distinct points, is called quasi-uniform if there is $c_u\in(0, 1]$ such that $c_u h_{\calx_n,\calx} \leq q_{\calx_n}$ for all $n\in\N$. 
\end{definition}
For quasi-uniform points there is a constant $c>0$ depending on $d$ and $\calx$ such that $h_{\calx_n, \calx}\leq c n^{-1/d}$.
We remark that it can be proven by a volume argument that $h_{\calx_n,\calx}\geq c' n^{-1/d}$ for any bounded set $\calx\subset \R^d$. It follows that quasi-uniform sequences of points provided an asymptotically optimal filling of $\calx$.

The approximation power of certain methods that we will discuss in the following will be measured by means of suitable sampling inequalities.
These formalize the idea that, for sufficiently regular functions, an approximant that is close to a target function on a sufficiently dense 
discrete set needs to also be close in the global sense \cite{Rieger17,Rieger20}. 

The results derived from such inequalities, which we recall in the following, differ according to the regularity of the kernel. 
To quantify this regularity, we will make use of Sobolev spaces $W_p^{\tau}(\calx)$ of possibly fractional smoothness $\tau\coloneqq k+s$, with $k\in\N$, $0<s< 
1$, and  
with $1\leq p\leq \infty$. 
Denoting as $D^\alpha$ the derivative with multiindex $\alpha\in\N_0^d$, if $s=0$ these spaces are equipped 
with the seminorm and norm defined as
\begin{equation*}
\seminorm{W_{p}^k(\calx)}{g}=\bigg(\sum_{|\alpha|=k}\lVert D^{\alpha}g\lVert^p_{L_p(\calx)}\bigg)^{1/p},\quad 
\norm{W_{p}^k(\calx)}{g}=\bigg(\sum_{|\alpha|\le k}\lVert D^{\alpha}g\lVert^p_{L_p(\calx)}\bigg)^{1/p},
\end{equation*}
while in the general case we set
\begin{equation*}
\norm{W_{p}^{k+s}(\calx)}{g}\coloneqq\left(\norm{W_p^k(\calx)}{g}^p +\seminorm{W_p^{k+s}(\calx)}{g}^p\right)^\frac1p,
\end{equation*}
where we used the seminorm
\begin{equation*}
\seminorm{W_p^{k+s}(\calx)}{g}
\coloneqq \left(\sum_{|\alpha|=k} \int_{\calx\times\calx} \frac{|D^\alpha g(x) -D^\alpha g(y)|^p }{\norm{2}{x-y}^{d+ps}}dxdy\right)^\frac1p,
\end{equation*}
with the usual modification for $p=\infty$.
We remark that if $\tau>d/2$ and $\calx$ is bounded and with a Lipschitz boundary, the Sobolev inequality gives the continuous embedding 
$W_2^{k+s}(\calx)\hookrightarrow C(\calx)$, and in particular all point evaluation functionals are continuous in $W_2^\tau(\calx)$, which is thus a RKHS 
according to Definition~\ref{def:rkhs}. 

We consider two classes of kernels of different smoothness, for both of which there are known sampling inequalities that we recall in the remaining part of 
this section.
First, if $K\in C^{k, s}(\calx)$ for some $k\in\N$, $s\in(0,1)$ with respect to both arguments, then there is a continuous embedding 
$\calh_K(\calx)\hookrightarrow C^{k,s}(\calx)$ (see 
e.g.~\cite{Wendland05}). Since $C^{k,s}(\calx)\subset W_2^{\tau}(\calx)$ if $\calx\subset\R^d$ is 
bounded and sufficiently regular, it also holds that 
$\calh_K(\calx)\hookrightarrow W_2^{\tau}(\calx)$.
Examples in this class are the Mat\'ern kernels and the compactly-supported Wendland kernels, see e.g. Chapter 10 \cite{Wendland05}. 
On the other hand, there are kernels of infinite smoothness whose RKHS is embedded in $W_p^k(\calx)$ for all $k\in\N$. Examples are the Gaussian and 
Inverse Multi Quadrics (IMQ) (see~\cite{Rieger10}).


The first result is Theorem 4.1 in~\cite{Arcangeli2007}, which generalized Theorem 2.6 \cite{WendlandRieger2005}, and it encompasses the case of Sobolev spaces. 
We report here the version of this result that is 
needed for our purposes.
We remark that we will use the result for $q=2, \infty$. In the second case the result of \cite{WendlandRieger2005} would be sufficient in our bounds, while we need this extended version for $q=2$.

\begin{theorem}[\cite{Arcangeli2007}]\label{thm:sampling_sobolev}
Let $k\in\N$, $s>0$, $1\le q \le \infty$, and $\tau\coloneqq k+s>d/2$, and let $m\in\N_0$ with $m \leq \lceil m_0\rceil -1$ and $m_0 = k+s - d (1/2 - 
1/q)_+$, where $(x)_+\coloneqq \max(x, 0)$.
Suppose that $\mathcal{X}\subset\mathbb{R}^d$ is bounded and satisfies an interior cone condition. 
Then there are constants $C, h_0>0$ such that for all $\calx_n\subset \calx$ with $h_{\calx_n,\calx}<h_0$, and any $z\in W_2^{\tau}(\calx)$ it holds 
\begin{equation*}
\seminorm{W_q^m(\calx)}{z}
\leq C \left(h_{\calx_n, \calx}^{\tau -m- d\left(\frac12 -\frac1q\right)_+}\seminorm{W_2^{\tau}(\calx)}{z} + h_{\calx_n, 
\calx}^{\frac d \gamma-m}\norm{2}{z_{|\calx_n}} 
\right),
\end{equation*}
where $\gamma\coloneqq \max(2, q)$, and $x/\infty=0$. 
\end{theorem}
\begin{remark}
A similar result holds if $\R^d$ is replaced by the sphere $\mathbb S^d$ (see Theorem 3.3 in~\cite{LeGia2006}), provided the fill distance is computed using the geodesic distance. We omit this result here for simplicity.
\end{remark}
For the second class of kernels, which are embedded in any Sobolev space, we will need instead the following. The result for Gaussians and IMQs was 
proved in Corollary 
5.1 \cite{Rieger10}, while Theorem 3.5 and Theorem 3.7 in~\cite{Lee2014a} extended the results to Multiquadrics and shifted surface splines. 

\begin{theorem}[\cite{Rieger10},\cite{Lee2014a}]\label{th:sampling_smooth}
Let $1\le q \le \infty$, and assume that $\mathcal{X}\subset\mathbb{R}^d$ is the union of finitely many compact cubes with equal side lengths. Let $\alpha\in\N_0^d$.
Then there are constants $c, c', h_0>0$ such that for all $\calx_n\subset \calx$ with $h_{\calx_n,\calx}<h_0$, and any $z\in \calh_K(\calx)$ the following 
holds:
\begin{itemize}
\item If $K$ is the Gaussian kernel, then
\begin{equation*}
\norm{L_q(\calx)}{D^\alpha z}\leq h^{-3|\alpha|}e^{c \frac {\log(h)} h} \norm{\calh_K(\calx)}{z} + h^{-2|\alpha|}e^{\frac {c'} h}\norm{\infty}{z_{|\calx_n}}.
\end{equation*}
\item If $K$ is the IMQ, Multiquadric, or shifted surface spline, then
\begin{equation*}
\norm{L_q(\calx)}{D^\alpha z}\leq h^{-3|\alpha|}e^{-\frac c h} \norm{\calh_K(\calx)}{z} +  h^{-2|\alpha|}e^{\frac {c'} h}\norm{\infty}{z_{|\calx_n}}.
\end{equation*}
\end{itemize}
%
\end{theorem}
The result is valid for much more general kernels, provided one can estimate the embedding constant of $\calh_K$ in $W_p^k$ for all $k\in\N$ and some $1\leq p< \infty$ (see Theorem 4.5 in \cite{Rieger10}). We stick here to this formulation for simplicity.

We remark that sampling inequalities have already been used in \cite{Krebs2009} to provide bounds on inverse problems, even if in a simpler setting when source conditions are considered.
Namely, \cite{Krebs2009} considers only Tikhonov regularization, pointwise bounded noise, and assumes that the problem is \emph{mildly-ill posed}, meaning that strong bounds on the reconstruction of $g$ can directly be transferred to strong bounds on the reconstruction of $f$ (see Assumption 4.6 in \cite{Krebs2009}).

\ifdefined\UNIGE
\begin{proposition}
    The solution in equation \eqref{approx_probl_sampling} coincides with the interpolant defined in equation \eqref{eq:interpolant} if $R_{\mathcal{Z}_n}(g')$ 
consists of the square loss, i.e. $R_{\mathcal{Z}_n}(g')=\sum_{i=1}^n (y_i-g'(x_i))^2$. 
\end{proposition}
\begin{proof}
     The solution in equation \eqref{approx_probl_sampling} admits the following representation
     \begin{equation}
         \hat{g}^{(n)}_{R}(x)=\sum_{i=1}^n \hat{c}_i K(x,x_i),
     \end{equation}
     where $\hat{c}=(\hat{c}_1,\dots,\hat{c}_n)$ is the solution of the following
     \begin{equation}
         \hat{c} = \arg\min_{c\in\mathbb{R}^n} \Vert y - \mathsf{K}c\Vert_2^2.
     \end{equation}
     Therefore, using the optimal condition, 
     \begin{equation}
         \mathsf{K}^T({y}-\mathsf{K}\hat{c})=0 \Longleftrightarrow \hat{c} = \mathsf{K}^{-1}y,
     \end{equation}
     i.e., $\hat{{c}}$ is the solution of the linear system derived from interpolation conditions in \eqref{eq:system}.
\end{proof}
\fi

\ifdefined\UNIGE
In $\mathcal{H}_K$, the approximation capability of $\hat{g}^{(n)}$ can be assessed by means of error bounds that may exploit various aspects of the 
interpolation design. Letting $K_x({\mathcal{X}_n}):=(K_x(x_1),\dots,K_x(x_n))$, a well-known pointwise error bound involves the so-called \textit{power 
function}
\begin{equation*}
    P_{K,n}=\lVert K_x - K_x({\mathcal{X}_n})\cdot(\mathsf{K}^{-1} K_x({\mathcal{X}_n}))\lVert_{\mathcal{H}_{K}},
\end{equation*} 
and reads as follows Theorem 14.2, p. 117 \cite{Fasshauer07}
\begin{equation}\label{eq:the_error}
    |g(x)-\hat{g}^{(n)}(x)|\le P_{K,n}(x)\lVert g\lVert_{H_K},\quad g\in H_K, \quad x\in\mathcal{X}.
\end{equation}
\textbf{Se g e' con rumore posso sostituire P con la kriging variance? (non ricordo dove l'ho letto)}
In the estimate \eqref{eq:the_error}, we observe that the error is split into a term that that relies on the underlying function $g$ and another that only 
depends on the nodes and on the kernel. Such bound can be improved as Section 9.3 \cite{Fasshauer15}
\begin{equation}\label{eq:the_error_bis}
    |g(x)-\hat{g}^{(n)}(x)|\le P_{K,n}(x)\lVert g-\hat{g}^{(n)}\lVert_{H_K},
\end{equation}
and to 
\begin{equation}\label{eq:the_error_tris}
    |g(x)-\hat{g}^{(n)}(x)|\le \delta P_{K,n}(x)\lVert \hat{g}^{(n)}\lVert_{H_K},
\end{equation}
if we assume $\lVert g-\hat{g}^{(n)}\lVert_{H_K}\le \delta\lVert \hat{g}^{(n)}\lVert_{H_K}$, i.e., the approximation of $g$ is \textit{good enough}.

\fi

\section{Main results}\label{sec:main_results}

This section is devoted to the main theoretical results of the paper.
In the first part, we show that spectral regularization yields an approximate interpolation operator whose behavior can be controlled through sampling inequalities.
In the second part, we apply these results to inverse problems and establish weak error estimates expressed in terms of projections onto operator-dependent test functions, thereby avoiding the use of source conditions.

\subsection{Approximate interpolation via spectral regularization}
We can now use the sampling inequalities of Section~\ref{sec:sampling} to obtain a bound on the norm of the error $e_\lambda^{(n)}\coloneqq g-\hat g_{s_{\lambda}}^{(n)}$, where $g\in \Range(A)$ and $\hat g_{s_{\lambda}}^{(n)}$ is computed by spectral regularization, as described in Section~\ref{sec:kernels_in_inv_prob} on the noise-free data, i.e.
\begin{equation}\label{eq:spectral_reg_g_noise_free}
    \hat{g}^{(n)}_{s_{\lambda}} = A s_{\lambda}(A_n^*A_n)A_n^* y.
\end{equation}
For this we need bounds on the terms $\seminorm{W_2^\tau(\calx)}{e_\lambda^{(n)}}$ and $\norm{2}{(e_\lambda^{(n)})_{|\calx_n}}$ appearing in Theorem~\ref{thm:sampling_sobolev} and Theorem~\ref{th:sampling_smooth}.
We derive these bounds in Lemma~\ref{lemma:sampling_spectral}, where we prove that both quantities can be controlled by terms depending solely on $\lambda$ and on the properties of the filtering function, and more specifically on the constants $D, C_{1/2}$ and $q$ (see Definition~\ref{def:regularization_fun}).
By inserting these bounds in the corresponding sampling inequalities we obtain the following error bounds for spectral regularization. 
We refer to Section~\ref{sec:sampling_spectral_proofs} for the proofs of these results.
\begin{proposition}\label{prop:sampling_spectral_finite_smooth}
Let $q$, $\tau$, $d$, $m$, $h_0$ $C$, $\calx$, $\calx_n$ be as in Theorem~\ref{thm:sampling_sobolev}, with $h_{\calx_n, \calx}\leq h_0$.
Assume that $\calh_K(\calx)\hookrightarrow W_2^\tau(\calx)$. Let $g\in \calh_K(\calx)$, and let $\hat g_{s_\lambda}^{(n)}$ be computed by spectral regularization with constants $D, C_{1/2}>0$ as in Definition~\ref{def:regularization_fun}. 
Then
\begin{equation}
\label{eq:bound_Lq}
\seminorm{W_q^m(\calx)}{g-\hat g_{s_\lambda}^{(n)}}
\leq C' h_{\calx_n, \calx}^{\frac d \gamma - m} \left(h_{\calx_n, \calx}^{\tau - \frac d2} + \lambda^{\frac12}\right) \norm{\calh_K}{g},
\end{equation}
with $\gamma\coloneqq\max(2,q)$, $C'\coloneqq C \max\{C_e(D+1), C_{1/2}\}$ and $C_e>0$ as in~\eqref{eq:norms_embedding}.
\end{proposition}

In the following, we will use the inequality in \ref{eq:bound_Lq} in the case $q=2$ and $q=\infty$.

\begin{proposition}\label{prop:sampling_spectral_infinite_smooth}
Let $\calx, \calx_n, c, c', h_0$ be as in Theorem~\ref{th:sampling_smooth}, with $h_{\calx_n, \calx}\leq h_0$.
Let $g\in\calh_K(\calx)$ and let $\hat g_{s_\lambda}^{(n)}$ be computed by spectral regularization with constants $D, C_{1/2}>0$ as in Definition~\ref{def:regularization_fun}.
Then, setting $C'\coloneqq \max\{D+1,  C_{1/2}\}$ we have the following:
\begin{itemize}
\item If $K$ is the Gaussian kernel, then
\begin{equation*}
\norm{L_q(\calx)}{D^\alpha\left(g-\hat g_{s_\lambda}^{(n)}\right)}
\leq C'h^{-3|\alpha|}e^{\frac {c'} h} \left(e^{-\frac{c'-c\log(h)}{h}} + \lambda^{1/2}\right)\norm{\calh_K(\calx)}{g}.
\end{equation*}
\item If $K$ is the IMQ, Multiquadric, or shifted surface spline, then
\begin{equation*}
\norm{L_q(\calx)}{D^\alpha\left(g-\hat g_{s_\lambda}^{(n)}\right)}
\leq C' h^{-3|\alpha|}e^{\frac {c'} h} \left(e^{-\frac{c'+c} h}+ \lambda^{1/2}\right)\norm{\calh_K(\calx)}{g}.
\end{equation*}
\end{itemize}
\end{proposition}

We remark that Lemma~\ref{lemma:sampling_spectral} generalizes Proposition 3.1 in~\cite{WendlandRieger2005}, which holds for Tikhonov regularization only. 
It follows that in this case Proposition~\ref{prop:sampling_spectral_finite_smooth} coincides with Proposition 3.6 in~\cite{WendlandRieger2005}, while Proposition~\ref{prop:sampling_spectral_infinite_smooth} is related with Theorem 6.1 in~\cite{Rieger10}.

\subsection{Weak error estimates for inverse problems}

Let $A$ be the compact operator defined in \eqref{eq:A_via_phi} described by the map $\phi$.
We now state two bounds for the weak error of spectral regularization: a general estimate valid for any bounded operator $A$, followed by an improved rate under a trace-class assumption on $A^*A$.
\begin{proposition}
\label{pro:psi_star_dagger}
    Let $\hat{f}^{(\delta,n)}_{s_{\lambda}}$ be the spectral regularized solution defined in equation \eqref{eq:spectral_reg_f}. 
    When $A$ is uniformly bounded, for each $\psi\in \Domain((A^*)^{\dagger})$, we have
     \begin{equation}
        \mathbb{E} |\langle \psi, \hat{f}^{(\delta,n)}_{s_\lambda}-f\rangle|_{\mathcal{H}_1} \le C_\psi \left(\Vert\hat{g}^{(n)}_{s_\lambda}-g\Vert_{\mathcal{H}_2} + \sigma_{\max} c
\frac{\nu}{\sqrt{n}}\frac{E}{\lambda}\right),
    \end{equation}
     where $E>0$ is the constant appearing in the second property of spectral regularization in Definition \ref{def:regularization_fun}, $\sigma^2_{\max}$ is the largest eigenvalue of $A^*A$, $c$ is is the uniform boundedness constant, and $\hat g^{(n)}_{s_{\lambda}}$ is the noise-free spectrally regularized solution introduced in \eqref{eq:spectral_reg_g_noise_free}.


In addition, if the operator $A^*A$ is of trace class, then the following bound holds:
    \begin{equation}
        \mathbb{E} |\langle \psi, \hat{f}^{(\delta,n)}_{s_\lambda}-f\rangle|_{\mathcal{H}_1} \le C_\psi \left(\Vert\hat{g}^{(n)}_{s_\lambda}-g\Vert_{\mathcal{H}_2} + \sigma_{\max}
\frac{\nu}{n}\frac{E}{\lambda} C''\right),
    \end{equation}
    where $C''>0$ is the constant bounding the trace norm of the operator $A_n^*A_n$.
\end{proposition}

We now derive analogous bounds for a different family of test functions commonly used in the study of sparse approximation~\cite{DeVore1996}, i.e.
\begin{equation*}
\mathcal{A}_1(\mathcal F)
\coloneqq\left\{\psi = \sum_{i\in \N} \alpha_i \phi_{x_i},~ x_i\in\calx, ~\sum_{i\in \N}|\alpha_i| < +\infty  \right\} \subseteq \mathcal F.
\end{equation*}

Since $\mathcal{A}_1(\mathcal F)$ is a subspace of the feature space $\mathcal F$, it is naturally endowed with a seminorm.
We extend this seminorm to $\mathcal{A}_1(\mathcal F) \oplus \Null(A)$ by setting $\seminorm{\mathcal A_1}{\psi+v}\coloneqq\seminorm{\mathcal A_1}{\psi}$ for any $v\in\Null(A)$, with $\psi\in \mathcal A_1(\mathcal F)$.
The statements below mirror the previous estimates, but now for test functions in $\mathcal A_1(\mathcal F)$, and are expressed in terms of the $\ell_1$-norm of the coefficient sequence $\alpha$ associated with a given $\psi \in A_1(\mathcal F)$.
\begin{proposition}
\label{pro:psi_l1}
Under the same hypothesis of proposition \ref{pro:psi_star_dagger}, we have
\begin{equation}
    \mathbb{E}(|\langle \psi, \hat{f}^{(\delta,n)}_{s_{\lambda}}-f\rangle|_{\mathcal{H}_1})\le
\Vert\alpha^\psi\Vert_{\ell_1}
\left(\Vert\hat{g}^{(n)}_{s_\lambda}-g\Vert_{\infty} + c^2 \frac{\nu}{\sqrt{n}}\frac{E}{\lambda} \right), 
    \end{equation}
for each $\psi\in \mathcal A_1(\mathcal F)$ whose coefficient sequence  $\alpha^\psi$ is absolutely summable.
Morevoer, if the operator $A^*A$ is of trace class then the following bound holds:
\begin{equation}
    \mathbb{E}(|\langle \psi, \hat{f}^{(\delta,n)}_{s_\lambda}-f\rangle|_{\mathcal{H}_1})\le
\Vert\alpha^\psi\Vert_{\ell_1}
\left(\Vert\hat{g}^{(n)}_{s_\lambda}-g\Vert_{\infty} + c \frac{\nu}{n}\frac{E}{\lambda} C''\right), 
    \end{equation}
\end{proposition}

By plugging the sampling inequalities into the propositions above, we obtain weak error bounds for inverse problems that are independent of any source conditions and controlled via suitable test functions.

\begin{proposition}(Weak error bounds for spectral regularization).
\label{pro:rates_psi_star_dagger}
   Let $\hat{f}^{(\delta,n)}_{s_{\lambda}}$ be the spectral regularized solution defined in equation \eqref{eq:spectral_reg_f} and $\mathcal{H}=L^2(\mathcal{X})$. For each $\psi\in \Domain((A^*)^{\dagger})$, we have the following weak error bound 
    \begin{equation}
        \mathbb{E}(|\langle \psi, \hat{f}^{(\delta,n)}_{s_\lambda}-f\rangle|_{\mathcal{H}_1})\le C_\psi \left(C_f h_{\calx_n, \calx}^{\frac d 2} \left(h_{\calx_n, \calx}^{\tau - \frac d2} + \lambda^{\frac12}\right) 
        + \sigma_{\max} c
\frac{\nu}{\sqrt{n}}\frac{E}{\lambda}\right),
    \end{equation}
    where $C_f = C' \Vert f^{\dagger}\Vert_{\mathcal{H}_1}$.
    Hence, the error becomes minimal for
    \begin{equation}
        \lambda \simeq \left(\frac{\nu^2}{n}\right)^{\frac{1}{3}} h_{\mathcal{X}_n,\mathcal{X}}^{-\frac{d}{3}}.
    \end{equation}
    If $A^*A$ is of trace class then
    \begin{equation}
        \mathbb{E}(|\langle \psi, \hat{f}^{(\delta,n)}_{s_\lambda}-f\rangle|_{\mathcal{H}_1})\le C_\psi \left(C_f h_{\calx_n, \calx}^{\frac d 2} \left(h_{\calx_n, \calx}^{\tau - \frac d2} + \lambda^{\frac12}\right) 
        + \sigma_{\max} 
\frac{\nu}{n}\frac{E}{\lambda} C''\right),
    \end{equation}
    where $C_f = C' \Vert f^{\dagger}\Vert_{\mathcal{H}_1}$.
    Hence, the error becomes minimal for
    \begin{equation}
        \lambda \simeq \left(\frac{\nu}{n}\right)^{\frac{2}{3}} h_{\mathcal{X}_n,\mathcal{X}}^{-\frac{d}{3}}.
    \end{equation}
\end{proposition}

The following corollary quantifies the optimal convergence rate of the weak error.

\begin{corollary}\label{cor:rates_psi_star_dagger}
   For each $\psi\in \Domain((A^*)^{\dagger})$, by assuming $h = O(n^{-\frac1d})$ we have the following optimal weak error bound 
 \begin{equation}
        \mathbb{E}(|\langle \psi, \hat{f}^{(\delta,n)}_{s_\lambda}-f\rangle|_{\mathcal{H}_1}) = O\left(\left(\frac{1}{n}\right)^{\min(\frac{\tau}{d}, \frac{1}{2})}\right), 
    \end{equation}
    with
    \begin{equation}
        \lambda = O(1)
    \end{equation}
If $A^*A$ is of trace class then we have
    \begin{equation}
        \mathbb{E}(|\langle \psi, \hat{f}^{(\delta,n)}_{s_\lambda}-f\rangle|_{\mathcal{H}_1}) = O\left(\left(\frac{1}{n}\right)^{\min(\frac{\tau}{d}, \frac{2}{3})}\right), 
    \end{equation}
    with
    \begin{equation}
        \lambda = O(n^{-\frac{1}{3}})
    \end{equation}
\end{corollary}

As in the previous case, we exploit the sampling inequalities to derive a weak error bound for test functions in $\mathcal A_1(\mathcal F)$.

\begin{proposition}(Weak error bounds for spectral regularization).
\label{pro:rates_psi_l1}
   Let $\hat{f}^{(\delta,n)}_{s_{\lambda}}$ be the spectral regularized solution defined in equation \eqref{eq:spectral_reg_f}. For each $\psi\in \mathcal A_1(\mathcal F)$
\begin{equation}
        \mathbb{E}(|\langle \psi, \hat{f}^{(\delta,n)}_{s_\lambda}-f\rangle|_{\mathcal{H}_1})\le \Vert\alpha\Vert_{\ell_1} 
\left(C_f \left(h_{\calx_n, \calx}^{\tau - \frac d2} + \lambda^{\frac12}\right) + c^2 \frac{\nu}{\sqrt{n}}\frac{E}{\lambda} \right), 
    \end{equation}
    where $C_f = C' \norm{\calh_K}{g}$.
    Hence the error becomes minimal for
    \begin{equation}
        \lambda \simeq \left(\frac{\nu^2}{n}\right)^{\frac{1}{3}}.
    \end{equation}
   If $A^*A$ is of trace class then
\begin{equation}
        \mathbb{E}(|\langle \psi, \hat{f}^{(\delta,n)}_{s_\lambda}-f\rangle|_{\mathcal{H}_1})\le \Vert\alpha\Vert_{\ell_1} 
\left(C_f \left(h_{\calx_n, \calx}^{\tau - \frac d2} + \lambda^{\frac12}\right) + c \frac{\nu}{n}\frac{E}{\lambda} C''\right), 
    \end{equation}
    Hence the error becomes minimal for
    \begin{equation}
        \lambda \simeq \left(\frac{\nu}{n}\right)^{\frac{2}{3}}.
    \end{equation}
\end{proposition}

The following corollary provides the optimal convergence rate in this case.

\begin{corollary}.
\label{cor:rates_psi_A1F}
  For each $\psi\in \mathcal A_1(\mathcal F)$, by assuming $h = O(n^{-\frac1d})$ we have the following optimal weak error bound 
 \begin{equation}
        \mathbb{E}(|\langle \psi, \hat{f}^{(\delta,n)}_{s_\lambda}-f\rangle_{\mathcal{H}_1}|) = O\left(\left(\frac{1}{n}\right)^{\min(\frac{\tau}{d}-\frac{1}{2}, \frac{1}{6})}\right), 
    \end{equation}
    with
    \begin{equation}
        \lambda = O(n^{-\frac{1}{3}}).
    \end{equation}
  If $A^*A$ is of trace class then
    \begin{equation}
        \mathbb{E}(|\langle \psi, \hat{f}^{(\delta,n)}_{s_\lambda}-f\rangle_{\mathcal{H}_1}|) = O\left(\left(\frac{1}{n}\right)^{\min(\frac{\tau}{d}-\frac{1}{2}, \frac{1}{3})}\right), 
    \end{equation}
    with
    \begin{equation}
        \lambda = O(n^{-\frac{2}{3}}).
    \end{equation}
\end{corollary}

Rates are the result of the well-known trade-off between bias and variance.
For instance, in Corollary \ref{cor:rates_psi_star_dagger}, which yields the strongest rates, the error is bias-dominated when $\frac{\tau}{d} < \frac{1}{2}$, so that increasing the kernel regularity $\tau$ improves the convergence rate. Once the regularity exceeds the threshold
$\frac{\tau}{d} \geq \frac{1}{2}$, however, the rate becomes variance-limited, and further increases in $\tau$ no longer improve it, as the error is now controlled by the amplification term. This behavior highlights a fundamental limitation: even with highly regular kernels, the convergence rate cannot exceed $O(n^{-1/2})$, reflecting the minimax optimal rate common to many statistical estimation problems. Consequently, while higher regularity is advantageous in low-dimensional settings, its benefit diminishes in higher dimensions or for fixed sample sizes.
This is the case, for instance, of the infinitely smooth Gaussian kernel: although it can, in theory, induce a very fast decay of the bias term, the overall convergence rate is ultimately limited by the variance. In high-dimensional settings or with a fixed sample size, this variance dominates the error, so that the high regularity of the  kernel does not necessarily translate into better practical performance.

\section{Applications}\label{sec:applications}
In this section we present several examples of linear functionals to which the error bounds shown above apply. 
In particular, we examine two classical cases: the inverse differentiation problem and the deconvolution problem with a translation-invariant kernel.

\subsection{Limited regularity}
Let $f \in \mathcal{H}_1 \coloneqq L^2[0,1]$. 
Consider the forward operator
\begin{equation}
 \label{eq:integral}   
g(x) = (Af)(x) =  \langle f, \phi_x \rangle_{\mathcal{H}_1} =
\int_{0}^x f(t) dt 
\end{equation}
where $\phi_x(t) = 1_{[0,x]}(t)$. 
The inverse problem is to find $f$ such that $g=Af$, given $g \in \mathcal{H}_2 \coloneqq L^2[0,1]$.
The image $R(A)$ is the Beppo Levi space, or homogeneous Sobolev space 
\cite{deny1954espaces}, equipped with the norm
$$
\Vert g\Vert_{\mathcal{H}_K} \coloneqq \| A^\dagger g \|_{\mathcal{L}_2} = \seminorm{W_2^1}{g} ~,
$$
where the generalized inverse $A^\dagger$ of the integral operator \eqref{eq:integral} is exactly the derivative $ D=A^\dagger$ over $R(A)$. Indeed, any function $g \in \mathcal H_2 \backslash  R(A)$ is such that $g(0) \neq 0$ and can be  approximated in the $\mathcal H_2$ norm with a sequence of elements in $R(A)$,  which means that $g \in \partial \mathcal R(A)$ and $R(A)^\perp = \emptyset$.
Moreover, as $g(0)=0$ for each $g\in R (A)$ no constant functions except zero are allowed in the range $R(A)$ and $A^\dagger = A^{-1}$.
Then, $R(A)$ is an Hilbert space with the Brownian reproducing kernel
\begin{equation}
    K(x_1,x_2) = \langle \phi_{x_1}, \phi_{x_2}\rangle_{L^2[0,1]} = \int_0^1 1_{[0,x_1]}(t) ~ 1_{[0,x_2]}(t) dt = \min(x_1,x_2)
\end{equation}
and scalar product
$$
\langle g_1, g_2 \rangle_{\mathcal H_K} \coloneqq \langle Dg_1, Dg_2 \rangle_{\mathcal{L}_2} ~.
$$

We now present two examples of test functions to which the weak error bounds  of Proposition 3 and Proposition 4 apply, respectively.
The first test function are the characteristic function $\chi_{xy} \coloneqq \phi_y - \phi_x = 1_{[x,y]}$ where $y>x$, we have a $| \chi_{xy} |_{\mathcal A_1} = 2$.
For these functions the bound in Proposition \ref{pro:psi_l1} readily applies, and the striking fact is that this bound is independent of the distance between $x$ and $y$, showing that we can localize the error with asymptotically infinite accuracy.

The second test functions 
are weakly differentiable functions on $[0,1]$ and belong to the domain of $(A^*)^\dagger$, since $A$ is the integral transform \eqref{eq:integral}.
In the following we report an example of  such a test function satisfying the condition of Proposition \ref{pro:psi_star_dagger}.
For a given $\varepsilon>0$, let $a,b\in [\varepsilon, 1-\varepsilon]$, $a< b$, and define 
\begin{equation*}
\psi_0\coloneqq \frac1{\varepsilon}\left(- \chi_{[a-\varepsilon, a]} + \chi_{[b, b+\varepsilon]}\right),
\end{equation*}
where $\chi$ is the indicator function. Since $\psi_0\in L_2(0,1)$, the function $\psi\coloneqq A^* \psi_0$ is in $R(A^*)$ 
by 
definition, and it holds
\begin{align*}
\psi(t) 
&= A^* \psi_0(t)
= \frac1\varepsilon\left(-\int_t^1\chi_{[a-\varepsilon, a]}(x)dx +\int_t^1 \chi_{[b, b+\varepsilon]}(x)  dx\right)\\
&= 
\frac1\varepsilon
\left(-\begin{cases}
\varepsilon, & t\leq a-\varepsilon,\\
a-t, & a-\varepsilon<t\leq a\\
0,& t>a
\end{cases}
+
\begin{cases}
\varepsilon, & t\leq b,\\
b+\varepsilon-t, & b<t\leq b+\varepsilon\\
0,& t>b+\varepsilon
\end{cases}\right)\\
&=
\begin{cases}
0, & t\leq a-\varepsilon,\\
1 - (a-t)/\varepsilon, & a-\varepsilon<t\leq a\\
1,& a<t\leq b\\
1 + (b-t)/\varepsilon, & b<t\leq b+\varepsilon\\
0,& t>b+\varepsilon,
\end{cases}
\end{align*}
which is a $C^0$-smoothed indicator function of the interval $[a, b]$. Moreover, 
\begin{equation*}
\norm{L_2(0,1)}{\psi_0}^2 
= \frac1{\varepsilon^2} \left(\int_{a-\varepsilon}^a dx + \int_b^{b+\varepsilon} dx \right)
= \frac2\varepsilon.
\end{equation*}
According to Lemma \ref{lemma:bound_h} and Remark~\ref{rem:c}, we thus have that 
\begin{equation*}
\inner{L_2(0,1)}{\psi, A^\dagger g} \leq {\sqrt{\frac2\varepsilon}} \norm{L_2(0,1)}{g}\; \fa g\in L_2(0,1).
\end{equation*}
This implies in particular that one can get arbitrarily close to measuring the error on a compactly supported area, up to accepting a constant growing as 
$\varepsilon^{-1/2}$.

\subsection{Infinitely smooth}

Now we see the case of a convolution operator with a infinitely smooth Gaussian kernel.
In this case, we show that bounds of Proposition \ref{pro:psi_l1} hold true for sufficiently regular approximation of any square integrable tezt function.
Suppose $A: \mathcal H_1 = \mathcal{L}^2(\mathbb{R}) \to \mathcal H_2 = \mathcal{L}^2(C)$ where $\mathcal C \subset \mathcal X$ a compact set, represents the space-invariant convolution operator
$$
Af (x) = \int_{\mathbb{R}} \tilde{\phi}(x-t) ~ f(t) ~ dt ~,
$$
where $\tilde{\phi}$ is a real-valued positive definite, continuous and square integrable function (cfr. equation 2.38 in \cite{bertero1998introduction}).
From the definition of the adjoint operator
$A^*g(t) = \int_C \phi(x+t)g(x) dx$,
we have that the reproducing kernel of $R (A)$ is 
$$
K(x,x') = \int_{\mathbb{R}} \tilde{\phi}(x-t) ~ \tilde{\phi}(x'-t) ~ dt
$$
is the autocorrelation function of $\tilde{\phi}$.
If we take 
$$
\phi_x(t) \coloneqq \tilde{\phi}(x-t)\coloneqq \sqrt{2/\pi} \exp(-(x-t)^2)
$$
the kernel $K$ is given by
$$
K(x,x') 
=
\sqrt{\frac{\pi}{2}}
\exp \left( {-\frac{1}{2} (x-x')^2} \right) ~.
$$
In this case $\mathcal A_1$ is dense in $\mathcal F = \mathcal L_2(\mathbb R)$ as it contains Gaussian functions with given variance. This result can be found in \cite{calcaterra2008linear}.

\begin{lemma}
\label{ThmBumps}
For any $\psi \in L^{2}\left(  \mathbb{R}\right)  $ and any
$\epsilon>0$ there exists $t_0>0$ and $N_\epsilon\in\mathbb{N}$ and $\alpha_{n}\in\mathbb{R}$ for $0<n \leq N_\epsilon$, such that for $t\neq 0$ and $|t|<t_0$
$$
\left\| \psi - 
\sum_{n=1}^{N_\epsilon} \alpha_{n} \exp{(-\left( ~\cdot~ -n t \right)^{2}) } \right\|_2 \leq \epsilon 
$$
\end{lemma}

Thus, for each $\psi \in \mathcal L_2(\mathbb R)$ and $\epsilon>0$, there exists a finite set of coefficient $\alpha_i$ with $i=1,\ldots,n$ so that the function
$$
\psi_{n_\epsilon} = \sum_{i=1}^{n_\epsilon} \alpha_i \phi_{x_i} 
$$
satisfies $
\left\| 
\psi_{n_\epsilon} - \psi 
\right\| \leq \epsilon$.
Therefore, we can approximate any functional $\psi$ with an arbitrary precision $\epsilon$ by $\psi_{n_\epsilon}$ and for such an approximation $\phi_\epsilon$, we have the error bound in proposition \ref{pro:psi_l1}.
The coefficients $\alpha_n$ depends also on the $t$ chosen to guaranteed the approximation, while $t$ can be chosen arbitrarily small, so that the approximation is valid even for compact subset $\mathcal C$.

\section{Proof of the approximation bound}\label{sec:proof_approx_err}

We start by proving a technical result to bound the value of the evaluation of a linear functional on $A^\dagger g$ in terms of the $\calh_2$-norm of $g\in 
\Domain(A^\dagger)$.

\subsection{Bounds on the $\calh_2$ norm}

We recall that $\Domain(A^\dagger)=\Range(A)\otimes\Range(A)^\perp$, $\Domain((A^*)^\dagger)=\Range(A^*)\otimes\Range(A^*)^\perp$.
\begin{lemma}\label{lemma:bound_h}
Let $\psi\in\calh_1$. 
There exists $C>0$ such that 
\begin{equation}\label{eq:bound}
\left|\inner{\calh_1}{\psi, A^\dagger g}\right|\leq C \norm{\calh_2}{g} \;\fa g\in \Domain(A^\dagger),
\end{equation}
if and only if 
\begin{equation}\label{eq:psi_in_range}
\psi\in \Domain((A^*)^\dagger)=\Range(A^*)+\Range(A^*)^\perp,
\end{equation}
and in this case~\eqref{eq:bound} holds with $C=\norm{\calh_2}{(A^*)^\dagger \psi}$.

Moreover, the condition \eqref{eq:psi_in_range} holds for any $\psi\in \calh_1$ if and only if the range $\Range(A)$ is finite dimensional, and in this
case in~\eqref{eq:bound} we can take
\begin{equation*}
C = \norm{\calh_1\to\calh_2}{(A^*)^\dagger}\norm{\calh_1}{\psi}.
\end{equation*}
\end{lemma}
\begin{proof}
First, if $g\in \Range(A)^\perp$ then $A^\dagger g=0$ and the bound holds trivially, so we can assume $g\in \Range(A)$.
Moreover, if $\psi\in \Range(A^*)^\perp$ we have by definition that $(A^*)^\dagger \psi= 0$ and thus $\norm{\calh_2}{(A^*)^\dagger \psi}=0$, and furthermore
$\Range(A^*)^\perp\subset \Null(A)$, so that $\inner{\calh_1}{\psi, A^\dagger g} = 0$ for all $g\in \Range(A)$ since $\Range(A^\dagger) \subset \Null(A)^\perp$. Thus in this
case~\eqref{eq:bound} 
and~\eqref{eq:psi_in_range} are clearly equivalent.

For all other cases, we first write all statements explicitly in terms of the bases.
Let $\psi\coloneqq \sum_{j=1}^\infty a_j u_j\in\calh_1$, $g\coloneqq \sum_j b_j v_j\in \Range(A)$ be generic elements, which means that it holds
\begin{equation}\label{eq:tmp_one} 
\norm{\calh_1}{\psi}^2=\sum_j a_j^2 < \infty,\;\;
\norm{\calh_2}{g}^2 =\sum_j b_j^2\leq \norm{\calh_K}{g}^2=\sum_j b_j^2/\sigma_j^2<\infty,
\end{equation}
which implies also that $A^\dagger g = \sum_j (b_j/\sigma_j) u_j$ is well defined.
Moreover, recall from~\eqref{eq:range_a_adj} that $\psi\in \Range(A^*)$ if and only if
\begin{equation}\label{eq:psi_in_range_a}
\sum_j a_j^2/\sigma_j^2<\infty,
\end{equation}
and in this case 
\begin{equation}\label{eq:tmp_two} 
(A^*)^\dagger\psi = (A^\dagger)^*\psi = \sum_j (a_j/\sigma_j) u_j, \;\;  \norm{\calh_2}{(A^*)^\dagger \psi}^2 = \sum_j (a_j/\sigma_j)^2.
\end{equation}
Finally, we have that
\begin{equation}\label{eq:tmp}
\inner{\calh_1}{\psi, A^\dagger g}=\sum_j \frac{a_j b_j}{\sigma_j}.
\end{equation}

To prove now that~\eqref{eq:psi_in_range} implies~\eqref{eq:bound}, we can simply apply the Cauchy-Schwartz inequality to~\eqref{eq:tmp} and get
\begin{equation*}
\inner{\calh_1}{\psi, A^\dagger g}
\leq \left(\sum_j(a_j/\sigma_j)^2\right)^{1/2}\cdot \left(\sum_jb_j^2\right)^{1/2}
= \norm{\calh_2}{(A^*)^\dagger \psi} \norm{\calh_2}{g},
\end{equation*}
where we used~\eqref{eq:tmp_one} and~\eqref{eq:tmp_two} in the last equality.

For the opposite direction assume instead that~\eqref{eq:psi_in_range} does not hold, which, thanks to \eqref{eq:psi_in_range_a}, implies that
\begin{equation*}
C_n\coloneqq\sum_{j=1}^n a_j^2/\sigma_j^2\xrightarrow{n\to\infty}+\infty.
\end{equation*}
We can consider now the function $g\in\calh_K=\Range(A)$ dedined by
\begin{equation*}
b_j\coloneqq\begin{cases}
            a_j/\sigma_j, &j\leq n\\
            0, &j>n,
            \end{cases}
\end{equation*}
which is indeed in $\calh_K$ for any finite $n\in\N$ (see~\eqref{eq:tmp_one}). However, we have from~\eqref{eq:tmp_one} and~\eqref{eq:tmp} that
\begin{align*}
\inner{\calh_1}{\psi, A^\dagger g}=\sum_{j=1}^n (a_j/\sigma_j)^2=C_n,\;\;
\norm{\calh_2}{g} =\left(\sum_{j=1}^n (a_j/\sigma_j)^2\right)^{1/2}=\sqrt{C_n},
\end{align*}
and thus~\eqref{eq:bound} may hold only if $C\geq C_n$. Since $n$ is generic and $C_n\to\infty$, there can be no finite $C$ that works 
in~\eqref{eq:bound}.

Finally, $\Range(A^*)+\Range(A^*)^\perp=\calh_1$ if and only if the sum is direct, which holds if and only if $\Range(A^*)$ is closed. This in turn is equivalent to require
that $\Range(A)$ is closed, and since $A$ is compact this is possible if and only if $\Range(A)$ is finite. In this case $(A^*)^\dagger$ is a bounded operator, and thus
the estimate on $C$ follows.
\end{proof}

From this lemma it is immediate to obtain a bound on the error in the approximation of the source, given a $\calh_2$-bound on the error in the approximation of  the data.
\begin{proposition}\label{prop:error_h}
Let $g\in\Range(A)$ and $\hat g \in \{ \hat{g}^{(n)}, \hat{g}_{\lambda}^{(n)}, \hat{g}^{(n,\delta)}_{s_{\lambda}}\}$ (see \eqref{approx_probl_sampling}, \eqref{approx_probl_sampling_noisy}and \eqref{eq:g_ldn}). Define $f\coloneqq A^\dagger g$ and $\hat f\coloneqq A^\dagger 
\hat g$. 
Then for each $\psi\in \Domain((A^*)^{\dagger}$ it holds
\begin{equation*}
\left|\inner{\calh_1}{\psi, f - \hat f }\right|\leq C \norm{\calh_2}{g-\hat g} \;\fa g\in \Domain(A^\dagger),
\end{equation*}
with $C\coloneqq\norm{\calh_2}{(A^*)^\dagger \psi}$.
\end{proposition}
\begin{proof}
The result follows by applying Lemma~\ref{lemma:bound_h} to the function $g-\hat g\in \Range(A)\subset\Domain(A^\dagger)$.
\end{proof}
We conclude the section with some remarks on this result.
\begin{remark}[On the proof of Lemma~\ref{lemma:bound_h}]
The direct implication can be proved also more directly, by observing that the expression $(A^\dagger)^*\psi$ is well defined if and only if $\psi \in 
\Domain((A^\dagger)^*)=\Domain((A^*)^\dagger)$, and in this case
\begin{align*}
\inner{\calh_1}{\psi, A^\dagger g}
= \inner{\calh_2}{(A^\dagger)^*\psi, g}
\leq \norm{\calh_2}{(A^\dagger)^*\psi} \norm{\calh_2}{g}
= \norm{\calh_2}{(A^*)^\dagger\psi} \norm{\calh_2}{g},
\end{align*}
and thus indeed~\eqref{eq:bound} works with $C=\norm{\calh_2}{(A^*)^\dagger\psi}$.
\end{remark}
\begin{remark}[Computation of $C$]\label{rem:c}
When~\eqref{eq:psi_in_range} holds, either $\psi\in \Range(A^*)^\perp$ and then $C=0$ as shown in the proof, or $\psi \in \Range(A^*)$. In this case, unless $\psi=0$,
there exists $\psi_0\in \Null(A^*)^\perp$ such that $\psi = A^* \psi_0$, and thus $C= \norm{\calh_2}{(A^*)^\dagger \psi} = \norm{\calh_2}{\psi_0}$.
\end{remark}

%
%

\subsection{Bounds on the $L_\infty$ norm}

We now turn to analogous bounds expressed in the $L_\infty$ norm.

\begin{lemma}\label{lemma:bound_l_inf}
Let $\psi\in\calh_1$ and assume that there are indices $I\subset \N$, points $\left(x_i\right)_{i\in I}\subset\calx$ and coefficients $\left(c_i\right)_{i\in I}\subset\R$, with
$\sum_{i=1}^\infty|c_i|\eqqcolon C< \infty$ such that
\begin{equation}\label{eq:approx_psi}
\lim_{n\to\infty}\norm{\calh_1}{\psi-\sum_{i=1}^n c_i \phi_{x_i}} = 0.
\end{equation}
Then
\begin{equation*}
\left|\inner{\calh_1}{\psi, A^\dagger g}\right|
\leq C \norm{L_\infty(\calx)}{g}\;\;\fa g\in \Domain(A^\dagger).
\end{equation*}
\end{lemma}
\begin{proof}
Set $\psi_n\coloneqq \sum_{j=1}^n c_j \phi_{x_j}$. We use the fact that $A f(x) = \inner{\calh_1}{\phi_x, f}$ for all $x\in \calx$ and for all $f\in \calh_1$ 
(see~\eqref{eq:A_via_phi}), and that $AA^\dagger g = g$ for all $g\in\calh_2$ to write
\begin{align*}
\left|\inner{\calh_1}{\psi_n, A^\dagger g}\right|
&= \left|\sum_{j=1}^n c_j \inner{\calh_1}{\phi_{x_j}, A^\dagger g}\right|
= \left|\sum_{j=1}^n c_j A (A^\dagger g)(x_j)\right|
= \left|\sum_{j=1}^n c_j g(x_j)\right|\\
&\leq\sum_{j=1}^n \left|c_j\right| \left|g(x_j)\right|
\leq C \norm{L_\infty(\calx)}{g},
\end{align*}
Taking the limit for $n\to\infty$ in the left-hand side concludes the proof thanks to~\eqref{eq:approx_psi} and the continuity of the norm.
\end{proof}

As in the previous paragraph, the lemma implies a result relating a weak error on $A^\dagger g$ with an $L_\infty$-error on $g$.

\begin{proposition}
Let $g\in\Range(A)$ and $\hat g \in \{ \hat{g}^{(n)}, \hat{g}_{\lambda}^{(n)}, \hat{g}^{(n,\delta)}_{s_{\lambda}}\}$ (see \eqref{approx_probl_sampling}, \eqref{approx_probl_sampling_noisy} and \eqref{eq:g_ldn}). Define $f\coloneqq A^\dagger g$ and $\hat f\coloneqq A^\dagger 
\hat g$. 
Then for each $\psi\in \mathcal A_1(\mathcal F)+ \Null(A)$ it holds
\begin{equation*}
\left|\inner{\calh_1}{\psi, f - \hat f }\right|\leq \seminorm{\mathcal A_1}{\psi} \norm{L_\infty(\calx)}{g-\hat g} \;\fa g\in \Domain(A^\dagger).
\end{equation*}
\end{proposition}
\begin{proof}
First, observe that $\Null(A)=\overline{\Range(A^\dagger)}^\perp\supset \Range(A^\dagger)^\perp$, thus if  $\psi\in \Null(A)$ it holds $\inner{\calh_1}{\psi, f-\hat f}=0$. We can then assume $\psi\in \mathcal A_1(\mathcal F)$.

The result then follows by applying Lemma~\ref{lemma:bound_l_inf} to the function $g-\hat g\in \Range(A)\subset\Domain(A^\dagger)$.
\end{proof}
\begin{remark}
Observe that in Proposition~\ref{prop:error_h} we do not need to add a term $v\in \Null(A)$ since $\Null(A)=\overline{\Range(A^\dagger)}^\perp=\Null((A^\dagger)^*)\subset \Domain((A^\dagger)^*)=\Domain((A^*)^\dagger)$, and thus this case is already covered by~\eqref{eq:bound} with $C=0$.
\end{remark}

\subsection{Sampling inequalities for spectral regularization}\label{sec:sampling_spectral_proofs}
We first recall how to compute interpolants using kernels that are positive definite, but non necessarily strictly positive definite.
\begin{lemma}\label{lemma:interp_pd}
Let $\calx_n\subset\calx$, and set $V(\calx_n)\coloneqq \Sp\{K(\cdot, x):x\in \calx_n\}$. 
Let $P_{\calx_n}:\calh_K(\calx)\to V(\calx_n)$ be the orthogonal projection onto $V(\calx_n)$.
Then, for $g\in \calh_K(\calx)$ the function $P_{\calx_n} g$ interpolates $g$ at $\calx_n$, and it can be written as
\begin{equation}\label{eq:proj_as_expansion}
P_{\calx_n} g(x)\coloneqq \sum_{i=1}^n \alpha_i K(x, x_i), \;\;x\in\calx,
\end{equation}
with $\alpha\coloneqq \mathsf{K}^{\dagger} g_{|\calx_n}$. 
\end{lemma}
\begin{proof}
The space $V(\calx_n)$ is finite dimensional, hence closed. It follows that the orthogonal projection $P_{\calx_n}$ is well defined, and by definition it holds that $g-P_{\calx _n} g\in V(\calx_n)^\perp$, and in particular 
\begin{equation}\label{eq:interp_cond}
0 = \inner{\calh_K(\calx)}{g - P_{\calx_n} g, k(\cdot, x_j)}
=g(x_j) - P_{\calx_n} g(x_j) \;\;\fa j=1, \dots, n.
\end{equation}
This means that $v\in V(\calx_n)$ coincides with $P_{\calx_n} g$ if and only if it interpolates $g$ at $\calx_n$. 

Now, since $P_{\calx_n} g\in V(\calx_n)$ it can be written in the form~\eqref{eq:proj_as_expansion} for some coefficient vector 
$\alpha\coloneqq(\alpha_i)_{i=1}^n\in \R^n$. The interpolation conditions~\eqref{eq:interp_cond} mean that $\mathsf{K}\alpha = g_{|\calx_n}$, and in 
particular $g_{|\calx_n}\in\Range(\mathsf{K})$. This in turns implies that the coefficients can be chosen as  $\overline\alpha\coloneqq \mathsf{K}^\dagger g_{|\calx_n}$, since  
$\mathsf{K}
\mathsf{K}^\dagger$ is the identity on $\Range(\mathsf{K})$ and thus $\mathsf{K} \overline \alpha = \mathsf{K} \mathsf{K}^\dagger g_{|\calx_n}=g_{|\calx_n}$.
\end{proof}

We can now prove a generalization of Proposition 3.1 in~\cite{WendlandRieger2005}.
\begin{lemma}\label{lemma:sampling_spectral}
Let $g\in\calh_K(\calx)$, $\calx_n\subset\calx$, $\lambda>0$, $s_\lambda$ a regularization function. 
Let $\hat g_\lambda^{(n)}\coloneqq \sum_{j=1}^n (\alpha_\lambda)_i K(\cdot, x_i)$ with $\alpha_\lambda = s_\lambda(\mathsf{K}) g_{|\calx_n}$. 
Then, if $s_\lambda$ has a qualification $q\geq 1/2$, for any $g\in \calh_K(\calx)$ it holds
\begin{align*}
\norm{2}{(g-\hat g_\lambda)_{|\calx_n}}&\leq C_{1/2} \lambda^{1/2} \norm{\calh_K}{g},\\
\norm{\calh_K(\calx)}{g-\hat g_\lambda^{(n)}}&\leq (D+1) \norm{\calh_K}{g},
\end{align*}
where $D, C_{1/2}$ are the constants of Definition~\ref{def:regularization_fun}.
\end{lemma}
\begin{proof}
We set $b\coloneqq g_{|\calx_n}$ and denote as $\mathsf{K}=U\Sigma U^T$ an SVD/eigendecomposition of $\mathsf{K}$, with $\Sigma\coloneqq \diag(\sigma^{(n)}_1, \dots, 
\sigma^{(n)}_n)$, $\sigma^{(n)}_i\geq \sigma^{(n)}_{i+1}\geq 0$. Furthermore we let $n'\leq n$ be the largest index for which $\sigma^{(n)}_{n'}>0$. 
We first observe that the norm of $P_{\calx_n} g$ is
\begin{align*}
\norm{\calh_K}{P_{\calx_n} g} 
&= \sqrt{\alpha^T \mathsf{K} \alpha}
= \sqrt{b^T \mathsf{K}^{\dagger} \mathsf{K} \mathsf{K}^{\dagger} b}
= \norm{2}{\mathsf{K}^{1/2} \mathsf{K}^{\dagger} b}\\
&= \norm{2}{U \Sigma^{1/2} U^T U \Sigma^\dagger U^T b}
= \norm{2}{\Sigma^{1/2} \Sigma^\dagger U^T b},
\end{align*}
where $\Sigma^{1/2} \Sigma^\dagger$ is a diagonal matrix. Moreover, note that multiplying by $\Sigma^{1/2}$ gives
\begin{equation}\label{eq:almost_inverse}
(\Sigma^{1/2} \Sigma^{1/2} \Sigma^\dagger)_{ii}
=
\begin{cases}
 \sqrt{\sigma^{(n)}_i} \sqrt{\sigma^{(n)}_i} / \sigma^{(n)}_i = 1, & 1\leq i\leq n'\\
 0, & n'+1\leq i\leq n.
\end{cases}
\end{equation}
For the discrete term, Lemma~\ref{lemma:interp_pd} gives
\begin{align*}
\norm{2}{(g-\hat g_\lambda)_{|\calx_n}}
&=\norm{2}{(P_{\calx_n} g-\hat g_\lambda)_{|\calx_n}}\\
&=\norm{2}{\mathsf{K}\alpha-\mathsf{K} \alpha_\lambda}\\
&=\norm{2}{\mathsf{K}\left(\mathsf{K}^\dagger -s_\lambda(\mathsf{K})\right)b}\\
&=\norm{2}{U\Sigma U^TU\left(\Sigma^\dagger - s_\lambda(\Sigma)\right)U^Tb}\\
&=\norm{2}{\Sigma \left(\Sigma^\dagger - s_\lambda(\Sigma)\right)U^Tb}.
\end{align*}
Now observe that~\eqref{eq:almost_inverse} implies that
\begin{equation*}
\Sigma \left(\Sigma^\dagger - s_\lambda(\Sigma)\right)
= \Sigma \left(\Sigma^\dagger - s_\lambda(\Sigma)\right) \Sigma^{1/2} \Sigma^{1/2} \Sigma^\dagger,
\end{equation*}
so we can continue with 
\begin{align*}
\norm{2}{(g-\hat g_\lambda)_{|\calx_n}}
&=\norm{2}{\Sigma \left(\Sigma^\dagger - s_\lambda(\Sigma)\right) \Sigma^{1/2} \Sigma^{1/2} \Sigma^\dagger U^Tb}\\
&\leq \norm{2}{\Sigma \left(\Sigma^\dagger - s_\lambda(\Sigma)\right) \Sigma^{1/2}}\norm{2}{\Sigma^{1/2} \Sigma^\dagger U^Tb}\\
&= \norm{2}{\Sigma \left(\Sigma^\dagger - s_\lambda(\Sigma)\right) \Sigma^{1/2}}\norm{\calh_K}{P_{\calx_n} g}.
\end{align*}
The norm of the first term can be computed explicitly since the matrix is diagonal, giving
\begin{align*}
\norm{2}{\Sigma \left(\Sigma^\dagger - s_\lambda(\Sigma)\right) \Sigma^{1/2}}
&= \max\limits_{t\in \{\sigma^{(n)}_1, \dots, \sigma^{(n)}_{n'}\}}\left|t \left(\frac1t - s_\lambda(t)\right)t^{1/2}\right|\\
&= \max\limits_{t\in \{\sigma^{(n)}_1, \dots, \sigma^{(n)}_{n'}\}}\left|t^{1/2} \left(1 - t s_\lambda(t)\right)\right|\\
&\leq C_{1/2} \lambda^{1/2}\;\;\fa \lambda>0, 
\end{align*}
where we used property (3) of Definition~\ref{def:regularization_fun}, with $q\geq 1/2$, in the last step.
Thus
\begin{equation*}
\norm{2}{(g-\hat g_\lambda)_{|\calx_n}}
\leq C_{1/2} \lambda^{1/2} \norm{\calh_K}{P_{\calx_n} g}
\leq C_{1/2} \lambda^{1/2} \norm{\calh_K}{g},
\end{equation*}
since $P_{\calx_n}$ is an orthogonal projection.

For the second term we proceed similarly and compute
\begin{align*}
\norm{\calh_K(\calx)}{\hat g_\lambda^{(n)}}
&= \sqrt{\alpha_\lambda^T \mathsf{K} \alpha_\lambda}
= \sqrt{b^T s_\lambda(\mathsf{K})\mathsf{K} s_\lambda(\mathsf{K}) b}
= \norm{2}{\mathsf{K}^{1/2} s_\lambda(\mathsf{K}) b}\\
&
= \norm{2}{\Sigma^{1/2} s_\lambda(\Sigma)  U^T b}
= \norm{2}{\Sigma^{1/2} s_\lambda(\Sigma) \Sigma^{1/2} \Sigma^{1/2} \Sigma^\dagger U^T b}\\
&\leq \norm{2}{\Sigma^{1/2} s_\lambda(\Sigma) \Sigma^{1/2}}\norm{\calh_K(\calx)}{P_{\calx_n} g},
\end{align*}
and now as before
\begin{align*}
\norm{2}{\Sigma^{1/2} s_\lambda(\Sigma) \Sigma^{1/2}}
&= \max\limits_{t\in \{\sigma^{(n)}_1, \dots, \sigma^{(n)}_{n'}\}}\left|t s_\lambda(t)\right|\leq D\;\;\fa \lambda>0,
\end{align*}
using (1) in Definition~\ref{def:regularization_fun}. 
Finally,
\begin{align*}
\norm{\calh_K(\calx)}{g-\hat g_\lambda^{(n)}}
&\leq \norm{\calh_K(\calx)}{g} +\norm{\calh_K(\calx)}{\hat g_\lambda^{(n)}}
\leq \norm{\calh_K(\calx)}{g} +D \norm{\calh_K(\calx)}{P_{\calx_n} g}\\
&\leq (D+1) \norm{\calh_K}{g},
\end{align*}
concluding the proof.
\end{proof}
These two bounds allow us to use a sampling inequality to bound the error of the spectrally regularized approximant.
For the first case we assume that there is a continuous embedding $\calh_K(\calx)\hookrightarrow W_2^{\tau}(\calx)$, i.e., there is $C_e>0$ such that
\begin{equation}\label{eq:norms_embedding}
\seminorm{W_2^\tau(\calx)}{g}
\leq 
\norm{W_2^\tau(\calx)}{g}
\leq C_e \norm{\calh_K(\calx)}{g}\;\; \forall g\in \calh_K(\calx).
\end{equation}
We can now prove Proposition~\ref{prop:sampling_spectral_finite_smooth} and 
Proposition~\ref{prop:sampling_spectral_infinite_smooth}.



\begin{proof}[Proof of Proposition~\ref{prop:sampling_spectral_finite_smooth}]
We use the sampling inequality of Theorem~\ref{thm:sampling_sobolev}, settitng $h\coloneqq h_{\calx_n, \calx}$ to simplify the notation.
Taking $v\coloneqq g-\hat g_\lambda^{(n)}$ and using the norm inequality~\eqref{eq:norms_embedding} and Lemma~\ref{lemma:sampling_spectral} we obtain
\begin{align*}
\seminorm{W_q^m(\calx)}{g-\hat g_\lambda^{(n)}}
&
\leq C \left(h^{\tau - m - d\left(\frac12 -\frac1q\right)_+}C_e \norm{\calh_K(\calx)}{g-\hat g_\lambda^{(n)}} + h^{\frac d \gamma - m}\norm{2}{(g-\hat g_\lambda^{(n)})_{|\calx_n}} 
\right)\\
&\leq C \left(h^{\tau - m -d\left(\frac12 -\frac1q\right)_+}C_e (D+1)\norm{\calh_K}{g} + h^{\frac d 
\gamma-m} C_{1/2} \lambda^{1/2} \norm{\calh_K}{g} \right)\\
&\leq C' \left(h^{\tau - m - d\left(\frac12 -\frac1q\right)_+} + h^{\frac d \gamma - m}  \lambda^{1/2}\right) \norm{\calh_K}{g}\\
&\leq C' h^{\frac d\gamma - m} \left(h^{\tau  - d\left(\frac12 -\frac1q\right)_+-\frac d \gamma} +  \lambda^{1/2}\right) \norm{\calh_K}{g},
\end{align*}
with $C'$ as in the statement. 
Observing that $d\left(\frac12 -\frac1q\right)_++\frac d \gamma = \frac d 2$ gives the first statement. 
\end{proof}

\begin{proof}[Proof of Proposition~\ref{prop:sampling_spectral_infinite_smooth}]
Using Theorem~\ref{th:sampling_smooth}, and following the same approach as in the proof of Proposition~\ref{prop:sampling_spectral_finite_smooth}, we get for the first case
\begin{align*}
\norm{L_q(\calx)}{D^\alpha\left(g-\hat g_\lambda^{(n)}\right)}
&\leq h^{-3|\alpha|}e^{c \frac {\log(h)} h} \norm{\calh_K(\calx)}{g-\hat g_\lambda^{(n)}} + h^{-2|\alpha|} e^{\frac {c'} h}\norm{\infty}{(g-\hat g_\lambda^{(n)})_{|\calx_n}}\\
&\leq h^{-3|\alpha|} \left(e^{c \frac {\log(h)} h} (D+1)\norm{\calh_K(\calx)}{g} + C_{1/2}\lambda^{1/2}e^{\frac {c'} h}\norm{\calh_K(\calx)}{g}\right)\\
&\leq C'h^{-3|\alpha|}e^{\frac {c'} h} \left(e^{-\frac{c'-c\log(h)}{h}} + \lambda^{1/2}\right)\norm{\calh_K(\calx)}{g},
\end{align*}
with $C'\coloneqq \max\{D+1,  C_{1/2}\}$.
The second case is completely analogous, giving
\begin{align*}
\norm{L_q(\calx)}{D^\alpha\left(g-\hat g_\lambda^{(n)}\right)}
&\leq h^{-3|\alpha|}e^{-\frac c h} \norm{\calh_K(\calx)}{g-\hat g_\lambda^{(n)}} + h^{-2|\alpha|} e^{\frac {c'} h}\norm{\infty}{(g-\hat g_\lambda^{(n)})_{|\calx_n}}\\
&\leq h^{-3|\alpha|}e^{\frac {c'} h} \left((D+1)e^{-\frac {c'+c} h}+ C_{1/2}\lambda^{1/2}\right)\norm{\calh_K(\calx)}{g},
\end{align*}
and using the same $C'$ as before.
\end{proof}

\begin{proof}[Proof of proposition \ref{pro:psi_star_dagger}]
Given $\hat{g}^{(\delta,n)}_{\lambda}$ in equation \eqref{eq:g_ldn}, we want to estimate the error 
$$
g-\hat{g}^{(\delta,n)}_{\lambda}
$$ 
as a function of $\lambda$, $\delta$ and $n$.
We have
\begin{eqnarray}  \label{eq:split_err_noise}  
    g-\hat{g}^{(\delta,n)}_{\lambda}
    &=& g- A_\lambda y^\delta \nonumber \\
    &=& g- A_\lambda y - A_\lambda \delta \nonumber \\
    &=& g- \hat{g}^{(n)}_{\lambda} - A_\lambda \delta \nonumber
\end{eqnarray}

By applying Proposition 7 to $\hat g = \hat{g}_\lambda^{(\delta,n)}$ we get
\begin{equation*}
\left|\inner{\calh_1}{\psi, f - \hat{f}_\lambda^{(\delta,n)} }\right|
\leq C \norm{\calh_2}{g-\hat{g}^{(\delta,n)}_{\lambda}}
\leq
\norm{\calh_2}{g- \hat{g}^{(n)}_{\lambda}} + \norm{\calh_2}{A_\lambda \delta}.
\end{equation*}
If $\calh_2=L_2(\calx)$, we take the expected value, and bound the second term using Proposition \ref{prop:noise_infty_l2}. 
\end{proof}

\begin{proof}[Proof of proposition \ref{pro:rates_psi_l1}]
Following the same approach, we apply instead Proposition 8 and obtain
\begin{align*}
\left|\inner{\calh_1}{\psi, f - \hat{f}_\lambda^{(\delta,n)}}\right|
&\leq \seminorm{\mathcal A_1}{\psi} \norm{L_\infty(\calx)}{g-\hat{g}_\lambda^{(\delta,n)}}\\
&\leq
\seminorm{\mathcal A_1}{\psi} 
\left(
\norm{L_\infty(\calx)}{g- \hat{g}^{(n)}_{\lambda}} + \norm{L_\infty(\calx)}{A_\lambda \delta}
\right)
\end{align*}
and bound the second term using Proposition~\ref{prop:noise_infty_l2}.
\end{proof}

\section{Proof of the noise amplification bound}\label{sec:proof_noise}

\begin{proposition}\label{prop:noise_infty_l2}
    Let $\hat{f}^{(\delta,n)}_{s_{\lambda}}$ be the spectral regularized solution defined in equation \eqref{eq:spectral_reg_f} and $\hat{f}^{(n)}_{s_{\lambda}}$ be the noisy-free spectral regularized solution defined as 
    \begin{equation}\label{eq:f_noisy_free}
\hat{f}^{(n)}_{s_{\lambda}}\coloneqq s_{\lambda}(A_n^*A_n)A_n^*A_n f.
    \end{equation}
    Then we have the following bounds
     \begin{equation}
        \mathbb E |(A(\hat{f}^{(\delta,n)}_{s_{\lambda}} - \hat{f}^{(n)}_{s_{\lambda}})(x)| \le  c^2 \frac{E\nu}{\sqrt{n}\lambda} 
    \end{equation}
    and
        \begin{equation}
        \mathbb E \Vert A(\hat{f}^{(\delta,n)}_{s_{\lambda}} - \hat{f}^{(n)}_{s_{\lambda}})\Vert_{\mathcal{H}_2} \le  \sigma_{\max} \frac{E \nu c}{\sqrt{n}\lambda},
    \end{equation}
    where $E>0$ is the constant of the second property of the spectral regularization, $\sigma^2_{\max}$ is the maximum eigenvalue of $A^*A$ and $c>0$ is the constant of the uniformly bounded hypothesis.\\
    If $A^*A$ is of trace class then we have the following bounds
     \begin{equation}
        \mathbb E |(A(\hat{f}^{(\delta,n)}_{s_{\lambda}} - \hat{f}^{(n)}_{s_{\lambda}})(x)| \le  c \frac{\nu}{n} \frac{E}{\lambda} C'',
    \end{equation}
    and
        \begin{equation}
        \mathbb E \Vert A(\hat{f}^{(\delta,n)}_{s_{\lambda}} - \hat{f}^{(n)}_{s_{\lambda}})\Vert_{\mathcal{H}_2} \le  \sigma_{\max} \frac{\nu}{n} \frac{E}{\lambda} C'',
    \end{equation}
    where $C''>0$ is the constant which bounds the trace norm of $A_n^*A_n$.
\end{proposition}

\begin{proof}

Let $\{(\sigma_l^{(n)},v_l^{(n)}\,u_l^{(n)})\}_{l=1}^n$ be the singular value decomposition of $A_n^*A_n$, then
\begin{equation}
    (A_n f)_i = \langle f, \phi_{x_i}\rangle = \sum_{l=1}^{n} \sigma_l^{(n)} \langle f, v_{l}^{(n)}\rangle (u_l^{(n)})_i
\end{equation}
and 
\begin{equation}
    \phi_{x_i} = \sum_{l=1}^n \sigma_l^{(n)} v_l^{(n)} (u_l^{(n)})_i.
\end{equation}
Let 
\begin{equation}
\label{learning spectral estimator}
\eta = \mathbf{y} - A_n f,
\end{equation} then
\begin{eqnarray*}
    s_{\lambda}(A_n^*A_n) A_n^*\eta &=& \sum_{j=1}^n s_{\lambda}((\sigma^{(n)}_j)^2) \langle v_j^{(n)}, A_n^*\eta\rangle_{\mathcal{H}_1} v_j^{(n)} = \\
    &=& \sum_{j=1}^n s_{\lambda}((\sigma^{(n)}_j)^2) \langle v_j^{(n)}, \frac{1}{n} \sum_{l=1}^n \phi_{x_l}\eta_l\rangle_{\mathcal{H}_1} v_j^{(n)} \\
    & = & \sum_{j=1}^n s_{\lambda}((\sigma^{(n)}_j)^2) \langle v_j^{(n)}, \frac{1}{n} \sum_{l=1}^n \eta_l\sum_{k=1}^n \sigma_k^{(n)} (u_k^{(n)})_l v_k^{(n)}\rangle_{\mathcal{H}_1} v_j^{(n)} \\
    & = & \frac{1}{n} \sum_{j=1}^n s_{\lambda}((\sigma^{(n)}_j)^2) \sum_{l=1}^n \sigma_j^{(n)} (u_j^{(n)})_l \eta_l v_j^{(n)} \\
    & = & \frac{1}{n} \sum_{j=1}^n s_{\lambda}((\sigma^{(n)}_j)^2)  \sigma_j^{(n)} \hat{\eta}_j v_j^{(n)},
\end{eqnarray*}
where $\hat{\eta}_j = \sum_{l=1}^n (u_j^{(n)})_l \eta_l = \langle \eta, u_j^{(n)}\rangle_{\mathbb{R}^n}$ is a random variable with mean 0 and variance $\nu^2$.
Therefore,
\begin{eqnarray*}
    |(A_{\lambda} \eta)(x)| = |\langle \phi_x, s_{\lambda}(A_n^*A_n) A_n^*\eta\rangle_{\mathcal{H}_1}| \le \Vert \phi_x\Vert_{\mathcal{H}_1} \Vert s_{\lambda}(A_n^*A_n) A_n^*\eta \Vert_{\mathcal{H}_1}
\end{eqnarray*}
\begin{eqnarray*}
    \mathbb{E}(\Vert s_{\lambda}(A_n^*A_n)A_n^*\eta\Vert_{\mathcal{H}_1}) & = &  \sqrt{\mathbb{E}(\Vert s_{\lambda}(A_n^*A_n)A_n^*\eta\Vert_{\mathcal{H}_1}^2)} \\
    & = & \sqrt{\frac{1}{n^2} \sum_{j=1}^n s^2_{\lambda}((\sigma^{(n)}_j)^2)  (\sigma_j^{(n)})^2 \mathbb{E}(\hat{\eta}^2_j) \Vert v_j^{(n)}\Vert_{\mathcal{H}_1}^2} \\
    & \le & \frac{E \nu}{n\lambda}\sqrt{\sum_{j=1}^n (\sigma_j^{(n)})^2} \le  \frac{E c \nu}{\sqrt{n}\lambda},
\end{eqnarray*}
where $E$ is the constant of the second property of the spectral regularization \eqref{property2 s_lambda} and we use $\sum_{j=1}^{n} (\sigma_j^{(n)})^2 = \sum_{j=1}^n \Vert\phi_{x_i}\Vert^2_{\mathcal{H}_1}\le n c$. 
In the case that the operator $A^*A$ is of trace class then by assuming $C''> 0 $ the constant which bounds the trace norm of the operator $A_n^*A_n$ the bound becomes $\frac{E C''\nu}{n\lambda}$.
Therefore for the uniformly boundness hypothesis we obtain
\begin{eqnarray*}
\mathbb{E}|A_{\lambda}\eta(x)|\le c^2 \frac{E \nu}{\sqrt{n}\lambda}.
\end{eqnarray*}
Further, we have the following bound
\begin{eqnarray}
       \mathbb{E}(\Vert A_{\lambda} \eta\Vert_{\mathcal{H}_2}) & \le & \Vert A \Vert  \sqrt{\mathbb{E}(\Vert s_{\lambda}(A_n^*A_n)A_n^*\eta\Vert^2_{\mathcal{H}_1})} \\
       & \le & \sigma_{\max} \frac{\nu E c}{\sqrt{n} \lambda} 
\end{eqnarray}
where $\sigma^2_{\max}$ is the maximum eigenvalue of $A^*A$.
In the case of the hypothesis that $A^*A$ is of trace class we obtain the following two bounds
\begin{eqnarray*}
\mathbb{E}|A_{\lambda}\eta(x)|\le c \frac{E C'' \nu}{n\lambda},
\end{eqnarray*}
\begin{eqnarray}
       \mathbb{E}(\Vert A_{\lambda} \eta\Vert_{\mathcal{H}_2}) & \le & \Vert A \Vert  \sqrt{\mathbb{E}(\Vert s_{\lambda}(A_n^*A_n)A_n^*\eta\Vert^2_{\mathcal{H}_1})} \\
       & \le & \sigma_{\max} \frac{\nu E C''}{n \lambda} .
\end{eqnarray}

\end{proof}

\section{Conclusions}
 \label{sec:cocnlusions}

We have presented an approach to error estimation in spectral regularization of inverse problems based on sampling inequalities. By interpreting the reconstruction error in a weak sense, we show that convergence rates can be established independently of a priori smoothness assumptions on the unknown solution. This is achieved by exploiting the connection between inverse problems and kernel approximation in RKHS, and by extending sampling inequalities to spectral regularization.
Our main contributions include: (i) spectral regularized error bounds expressed in terms of the fill distance of sampling points for kernel approximation problems; (ii) point-wise variance bounds beyond traditional $L_2$ norm, in the framework of kernel approximation with additive Gaussian noise; and (iii) weak convergence estimates for spectral regularized solutions of inverse problems depending on the fill distance for two given classes of test functions according to the forward operator.
Our findings show that increasing kernel regularity improves convergence rates only up to a certain degree. For sufficiently smooth kernels or high dimensions, the rate becomes variance-limited and cannot exceed $O(n^{-\frac{1}{2}})$, yielding the classical minimax optimal rate for many statistical estimation problems.
Importantly, the bounds depend on geometric properties of the sampling set and are independent of the regularity of the unknown solution, thus allowing possibly relevant implications for  practical applications, such as image reconstruction and kernel-based learning, where sampling patterns are often known.


\bmhead{Acknowledgements} All authors acknowledge the Gruppo Nazionale per il Calcolo Scientifico - Istituto Nazionale di Alta Matematica (GNCS - INdAM). This research was conducted within the framework of the MIUR Excellence Department Project awarded to the Department of Mathematics, University of Genoa (CUP D33C23001110001).
G.S. was partially supported by the project ``Perturbation problems and asymptotics for elliptic differential 
equations: variational and potential theoretic methods'' funded by the program ``NextGenerationEU'' and by MUR-PRIN, grant 2022SENJZ3.

\bibliography{references}


\end{document}